\documentclass[12pt, oneside]{amsart}
\usepackage{amscd}
\usepackage{amssymb}
\usepackage{remark}
\usepackage{enumerate}

\input BoxedEPS.tex
\SetRokickiEPSFSpecial
\HideDisplacementBoxes

\newtheorem{thm}{Theorem}[section]

\newtheorem{lemma}[thm]{Lemma}
\newtheorem{lem}[thm]{Lemma}
\newtheorem{cor}[thm]{Corollary}

\newtheorem{prop}[thm]{Proposition}
\newtheorem{proposition}[thm]{Proposition}
\newremark{rem}[thm]{Remark}
\newremark{ques}[thm]{Question}
\newremark{remark}[thm]{Remark}
\newremark{exmp}[thm]{Example}
\newremark{example}[thm]{Example}
\newremark{defn}[thm]{Definition}

\newremark{emp}[thm]{\kern -4pt}

\newremark{fact}[thm]{Fact}
\newremark{Facts}[thm]{Facts}

\newcommand{\Spec}{\mbox{\rm Spec}\kern 1pt}

\newcommand{\ba}{\begin{array}}

\newcommand{\ea}{\end{array}}

\newcommand{\xxx}{{\mathcal X}}
\newcommand{\zzz}{{\mathcal Z}}
\newcommand{\yyy}{{\mathcal Y}}
\newcommand{\ooo}{{\mathcal O}}

\newcommand{\bQ}{{\mathbb Q}}
\newcommand{\Qb}{{\overline \bQ}}
\newcommand{\bZ}{{\mathbb Z}}

\newcommand{\Zp}{\bZ_p}

\newcommand{\cm}{{\mathcal M}}
\newcommand{\cO}{{\mathcal O}}

\newcommand{\Chab}{\textnormal {Chab}}

\newcommand{\wto}{\overline{\omega}}

\newcommand{\cU}{{\mathcal U}}
\newcommand{\cP}{{\mathcal P}}

\newcommand{\lra}{\longrightarrow}
\newcommand{\cC}{{\mathcal C}}
\newcommand{\cX}{{\mathcal X}}

\newcommand{\cXb}{\overline \cX}

\newcommand{\rank}{\textnormal{rank}}

\newcommand{\Proj}{\textnormal{Proj }}
\newcommand{\Qp}{\bQ_p}

\newcommand{\cS}{{\mathcal S}}

\newcommand{\Fq}{{\mathbb F}_q}

\begin{document}

\pagestyle{headings}

\title{ Thue Equations and the method of Chabauty-Coleman }
\author{Dino Lorenzini
 and Thomas J. Tucker}

\thanks{The authors wish to thank John M. Grigsby
for his support of research in number theory.  They would also like to
thank W. McCallum and J. Wetherell for helpful mathematical conversations.}


\date{\today}

\begin{abstract}

Let $F(x, y) \in {\mathbb Z}[x, y]$ be a homogenous polynomial
of degree $n$, assumed to have unit content.
  Let $h \in {\mathbb Z}$ be such that the
polynomial $hz^n - F(x, y)$ is irreducible in
$\overline{\mathbb Q}[x, y, z]$.  We denote by $X_{F,h}/{\mathbb Q}$ the
nonsingular complete model of the projective plane curve $C_{F,h}/{\mathbb Q}$
defined by the equation $hz^n - F(x, y) = 0$.  
A classical
Thue equation is an equation $F(x, y) = h$ where $F(x, 1)$
does not have repeated roots.  
  Let 
$N(F, h) $ denote the cardinality of the set
$$ \{(x,
y) \in {\mathbb Z}^2 \mid F(x, y) = h \ {\rm and \ } \ {\rm gcd}(x, y) =
1\}. $$ 
Let $r(X_{F,h})$ denote the Mordell-Weil rank over ${\mathbb Q}$ of
the jacobian of $X_{F,h}/{\mathbb Q}$.
Our main theorem is:

\medskip
\noindent {\bf Theorem. \ } {\it Assume that the genus $g$ of
$X_{F,h}$ is at least $2$. 
If $r(X_{F,h}) <  g$, then
$
N(F, h) \leq  2n^3-2n -3.$
}\medskip

To prove this theorem, we first refine
the method of Chabauty-Coleman so that it can be applied to any
regular model of a curve $X_{F,h}/{\mathbb Q}_p$ with $p>n$. (Coleman's
Theorem applies only to curves having good reduction at
a prime $p > 2g$.)  We then  describe some regular
open subsets of a normal model of the curve $X_{F,h}/{\mathbb Q}_p$ and  prove
our main theorem on primitive solutions of Thue equations using this
model. We also present some refinements of our main theorem
in some special cases
  where we obtain a bound of the form $N(F, h) \leq  O(n^2)$,
and discuss examples of twists of Fermat curves.
\end{abstract}

\maketitle

Let ${\mathcal O}_K$ be a domain with field of fractions $K$.
Let $F(x, y) \in {\mathcal O}_K[x, y]$ be a homogenous polynomial
of degree $n$, assumed to have unit content (i.e., the coefficients of $F$ generate
the unit ideal in ${\mathcal O}_K$). Assume that $\gcd(n, {\rm char}(K)) =1$.
 Let $h \in {\mathcal O}_K$ and assume that the
polynomial $hz^n - F(x, y)$ is irreducible in
$\overline{K}[x, y, z]$.  We denote by $X_{F,h}/K$ the
nonsingular complete model of the projective plane curve $C_{F,h}/K$
defined by the equation $hz^n - F(x, y) = 0$.  Recall that
$hz^n -F(x, y) = 0$ defines a nonsingular curve
if $F(x, 1)$ has simple roots in $\overline{K}$, and that if
$F(x, 1) = c \prod^s_{i=1}(x - \alpha_i)^{n_i}$ with $\prod_{i \neq j} (\alpha_i - \alpha_j) \neq 0$, then the genus $g$ of
$X_{F,h}$ is given by the formula
$$
2g - 2 = n(s - 2) - \sum^s_{i=1} \gcd(n, n_i).
$$
We shall assume in this article that $g \geq 2$.  When $K$ is a
number field, Mordell's Conjecture implies that $|X_{F,h}(K)| <
\infty$. Caporaso, Harris, and Mazur (\cite{CHM}, 1.1) have shown
that if Lang's conjecture for varieties of general type is true, then
 for any number field $K$, the size $|X(K)|$ of the set
of $K$-rational points of any curve $X/K$ of genus $g(X) \ge 2$ can be
bounded by a constant depending only on $g(X)$.
%% Changed the wording of the above sentence slightly
  Prior to the paper \cite{CHM}, Mazur
and others had asked whether $|X(K)|$ can be bounded by a constant
depending only on $g(X)$ and the Mordell-Weil rank of $X/K$ over $K$
(that is, the rank of the group $J(K)$ of $K$-rational points of the
jacobian $J/K$ of $X/K$).  These far-reaching questions are totally
open.  As we shall recall in section \ref{sec.Chabauty}, the method of
Chabauty-Coleman sometimes yields a bound for $|X_{F,h}(K)|$ depending only
on $g(X_{F,h})$ when it is known in advance that the Mordell-Weil rank of
$X_{F,h}/K$ is small.  Unfortunately, the Chabauty-Coleman method does not
yield a bound for $|X_{F,h}(K)|$ independent of the coefficients of
$hz^n - F(x, y)$ for all curves of the form $X_{F,h}$.  It does,
however, produce such a nice bound for the number of primitive
integral points of $X_{F,h}/K$, as we now explain.

Let $K = {\mathbb Q}$ and ${\mathcal O}_K = {\mathbb Z}$.  A classical
Thue equation is an equation $F(x, y) = h$ where $F(x, 1)$
does not have repeated roots.  Thue showed in 1909 that such
an equation has finitely many solutions $(x, y) \in {\mathbb
Z}^2$ if $n \ge 3$. Let us say that $(x,y)$ is a primitive solution
if $\gcd(x,y) =1$.
In this work, we are
interested in the following open question raised for instance 
by Erd\"os, Stewart, and Tijdeman (\cite{Ste}, p.\ 816).  Let 
$N(F, h) $ denote the cardinality of the set
$$ \{(x,
y) \in {\mathbb Z}^2 \mid F(x, y) = h \ {\rm and \ } \ {\rm gcd}(x, y) =
1\}. $$ Is there a bound for $N(F,
h)$ in terms of $n$ only whenever $n \ge 3$?  Two
known results on $N(F, h)$ are as follows:
 
\medskip
\noindent {\bf Theorem} {\rm (Bombieri-Schmidt, \cite{B-S}).}  
{\it Assume that $F(x,1)$ has distinct roots in $\overline{\mathbb
    Q}$.  There exists a constant $B_1$, which can be taken to be 215
  when $n$ is sufficiently large, such that $N(F, h) \le B_1 n^{w(h) +
    1}$, where $w(h)$ equals the number of prime factors of $h$.}
\medskip

This bound depends on $n$ and $h$.
Let $r(X_{F,h})$ denote the Mordell-Weil rank over ${\mathbb Q}$ of
the jacobian of $X_{F,h}/{\mathbb Q}$.

\medskip
\noindent
{\bf Theorem} {\rm   (Silverman, \cite{Sil}).  }{ \it 
Assume that $F(x,1)$ 
has distinct roots in $\overline{\mathbb Q}$.
There exists an 
ineffective
constant $h(F)$ such that for all $n^{th}$ power-free $h >
h(F)$, then
$$
N(F, h) \le n^{2n^2}(8n^3)^{r(X_{F,h})}.$$
} \medskip

For a fixed $F$, this bound depends only on $n$ and $r(X_{F,h})$,
but only works for $h$ sufficiently large. Our main theorem is:

\medskip
\noindent {\bf Theorem   \ref{thm.main} }{\it 
%Assume that $n \ge 4$, and that $F(x,1)$  
%has distinct roots in $\overline{\mathbb Q}$. 
If $r(X_{F,h}) <  g(X_{F,h})$, then
$
N(F, h) \leq  2n^3-2n -3.$
}\medskip

This bound only holds when $r(X_{F,h})$ is small, but when it
holds, it depends only on $n$. In section \ref{sec.ref}, we are able
to refine our method in some special cases
  to obtain a bound of the form $N(F, h) \leq  O(n^2).$
%We also consider cases where $F(x,1)$ has repeated roots.

%% Tom
The theorems of \cite{B-S} and \cite{Sil} both hold when $n=3$.
In this case, $X_{F,h}$ is an elliptic curve with potentially good reduction.
When its rank is assumed to be zero, that is, under the hypothesis of Theorem \ref{thm.main}, it is easy to find a bound for the order of the finite group $X_{F,h}({\mathbb Q}) $
 by reducing modulo $p$ for appropriate primes $p$.

Both \cite{B-S} and \cite{Sil} make use of the Thue-Siegel-Roth
theorem on approximations of algebraic numbers.  Typically, this
theorem yields a bound on the number of solutions to an equation that
depends on the size of the coefficients defining the equation.
Bombieri and Schmidt transform the coefficients of Thue equations with
elements of $\textnormal{SL}_2(\bZ)$ to overcome this difficulty in
the case $h=1$; the general case of their theorem is proved via
induction on the number of prime factors of $h$.  Silverman works with
a fixed $F$ and uses the Thue-Siegel-Roth theorem to dispose of
solutions to $F(x,y) = h$ corresponding to the so-called large
approximations of the roots of $F$.  The number of these large
approximations can only be bounded in terms of the size of the
coefficients of $F$, but, in any case, there are only finitely many
such approximations for a fixed $F$ and thus only finitely $h$ for
which these approximations contribute to the number of solutions to
$F(x,y) = h$.

The proof of Theorem \ref{thm.main}, by contrast, does not involve
diophantine approximation.  As mentioned earlier, Theorem
\ref{thm.main} is proved using the method of Chabauty-Coleman on the
curve $X:= X_{F,h}/{\mathbb Q}$.  In order to use this method to bound $|X({\mathbb Q})|$,
 one needs to pick
a prime and compute enough of a regular model
${\mathcal X}/{\mathbb Z}_p$ of $X/{\mathbb Q}_p$ to be able to bound the
number $N_1$ of components of multiplicity 1 in the special fiber
$\overline{\mathcal X}/{\mathbb F}_p$.  
The number $N_1$ is not, in general, bounded by a constant depending only on $g(X)$ 
(see \ref{rem.toomanyP1}).  Hence, this method does
not always enable us to bound $|X(\bQ)|$
in terms of $g(X)$ only\footnote{It may be possible to bound the number $N_1$
in terms of the size of $h$ and the coefficients of $F$. If such were the case,
then the method of Chabauty-Coleman would provide a bound for $|X(\bQ)|$
in terms of $n$ and the size of $h$ and the coefficients of $F$.}.
Surprisingly, however, it is possible to bound, in terms
of $n$ only, the
number of reduction classes 
in the special fiber of a  regular model ${\mathcal X}/\Zp$ 
of the  integer solutions $(x,y)$ of $F(x,y) = h$ with
$\gcd(x,y) =1$. 
Note that there are only few families 
of curves $X_n/{\mathbb Q}$ with lim sup $g(X_n)= \infty$  for which  a bound 
for $|X_n({\mathbb Q})|$ is known and depends only on $g(X_n)$.
For two such families,   modular curves (work of  Mazur) and   Fermat curves (work of Wiles), the sets $X_n(\mathbb Q)$ of rational 
points are in fact
completely described.

The Mordell-Weil rank of a  curve $X/K$ is in general very hard
to compute.  The case of superelliptic curves of the form
$y^p = f(x)$ with $p \mid{\rm deg}(f)$ is treated in \cite{P-S}.
An  effective algorithm in the case $p = 2$ and ${\rm deg}(f)=6$
has been
implemented with Magma \cite{Sto}, and  this algorithm can be used
 to produce some explicit examples 
with $n=6$ where the
bound given in Theorem \ref{thm.main} holds. Indeed, as is recalled in \ref{tobedone},
the method of Chabauty-Coleman can be applied to the curve $X_{F,h}$
with $n=6$
if its hyperelliptic quotient $hz^2 = F(x,1)$ has Mordell-Weil rank
at most $1$.

We do not discuss here the problem of effectively determining all of the
solutions to a Thue equation.  For some recent work on this problem,
we refer the reader to \cite{BH} and \cite{TdW}. 

This paper is organized as follows.  In the first section, we refine
the method of Chabauty-Coleman so that it can be applied to any
regular model of a curve $X_{F,h}/{\mathbb Q}_p$ with $p>n$. (In
\cite{Co2}, the method applies only to curves having good reduction at
a prime $p > 2g$.)  In the second section, we first describe some
regular models of the curves $X_{F,h}/{\mathbb Q}_p$ and we then prove
our main theorem on primitive solutions of Thue equations using these
models. In the last section, we present some refinements and discuss
examples of twists of Fermat curves.

\section{The method of Chabauty-Coleman} \label{sec.Chabauty}

We begin by fixing some notation for this section.  
Let $K$ be any number field with a place $v$ over a prime $p$.
Let $K_v$ denote the completion of $K$ at $v$, with uniformizer $\pi$
and residue field   $\Fq$, where $q$ is a power of $p$.
Let $X/K$ be
any nonsingular, connected, complete, curve of genus $g$. Let $\cX/{\mathcal O}_{K_v}$ be 
a proper   model for $X/K_v$. 
Denote by $\overline{\mathcal X}$
the special fiber 
${\mathcal X}\times_{\text{\rm Spec}
({\mathcal O}_{K_v})}\Spec({\mathcal O}_{K_v}/(\pi)) $.
Let $J/K$ be the Jacobian of $X$.   
Since $\cX$ is proper, we have a reduction map
$$ r:\cX(\overline{K}_v) \lra \cXb({\overline{\mathbb F}_q}),$$
which sends points in $\cX(K_v)$ to points in $\cXb(\Fq)$.
If $Q \in \cXb(\Fq)$,  denote by $D_Q$ the set $r^{-1}(Q)$.
When $P\in X(K_v)$, the set $D_{r(P)}$ is called the residue class of $P$.

Let $A/K$ be any abelian variety. 
As a $p$-adic Lie group, $A(K_v)$ is endowed with
a logarithm map $\text{\rm log}: A(K_v) \to T_0(A)(K_v)$,
where $T_0(A)(K_v)$ denote the tangent space to $A$ at $0$.
For an abelian variety $A$ over a number field $K$ with completion
$K_v$,  the {\it Chabauty rank of $A$ at $v$} is
the integer
$$\Chab(A,K,v) := \dim_{K_v} (\log(A(K))\otimes_{\mathbb Z} K_v).$$
Note that since $\log$ is a homomorphism, $\Chab(A,K,v)$ is less than
or equal to the Mordell-Weil rank of $A(K)$.

Let us  now recall briefly the main ideas of the method of Chabauty-Coleman.
Consider again $A(K_v)$ as a $p$-adic Lie group.
Given any
  global differential $\eta \in \Gamma(A,\Omega_{A/K_v})$, 
there is a unique $p$-adic analytic homomorphism 
$$\lambda_\eta : A(K_v) \lra K_v$$
such that $d(\lambda_\eta) = \eta$
(see \cite[Lemma 1.3.2]{Weth}).  In fact, it follows from standard
facts on formal groups (see \cite[Theorem 1]{Freije}), for example)
that one can take $\lambda_\eta$ to be
$$(\eta)_0 \circ \log: A(K_v) \lra K_v,$$
where $(\eta)_0$ is the stalk at $0$ of $\eta$ considered as a dual
element of the tangent bundle of $A$.  
The key  to the method
of Chabauty-Coleman is the remark that when $\Chab(A,K,v) < \text{\rm dim}(A)$,
then there exists a differential ${\eta}$ such that the 
homomorphism $\lambda_{\eta}: A(K_v) \to K_v$ vanishes on $A(K)$.

Suppose now that $A/K$ is the Jacobian $J/K$
of a curve $X/K$. Define the {\it Chabauty rank of $X/K$ at $v$}
to be $\Chab(J,K,v)$, and denote it by $\Chab(X,K,v)$.  Assume that
 there is an embedding $j:X \lra J$ defined over $K$.
This embedding induces an
isomorphism $j^*$ from $\Gamma(A,\Omega_{A/K_v})$ to
$\Gamma(X,\Omega_{X/K_v})$.  Since every $\omega \in
\Gamma(X,\Omega_{X/K_v})$ is $j^* (\eta)$ for some $\eta \in
\Gamma(A,\Omega_{A/K_v})$, there is an analytic map 
\begin{equation} \label{eq.lambda}
 \lambda_\omega: X(K_v) \lra K_v
\end{equation}
such that 
$d(\lambda_\omega) = \omega$.
Furthermore, for any $P \in \cX(K_v)$, there is a power series that
converges on the residue class of $P$ and is equal to $\lambda_\omega$
as a function on that residue class (see \cite[Lemma
1.7.3]{Weth}).
%%Tom  1.7.3 is correct? you had 1.7.4
 When $\cX/\cO_{K_v}$ is a regular model
for $X/K_v$  and $\omega$ is in
$\Gamma(\cX, \Omega_{\cX/\cO_{K_v}})$, 
%%Tom why do we need that?
%%(considered as a lattice in
%%$\Gamma(X, \Omega_{X/K_v})$, 
the power series has the form
$$ \lambda_\omega = a_0 + \sum\limits_{m=1}^{\infty} \frac{a_m}{m}
  u^m,$$
where all of the $a_i$ are in $\cO_{K_v}$, and $u: D_{r(P)} \to \pi {\mathcal O}_K$ is a local coordinate.  This power series is
obtained by formally integrating a power series expansion for
$\omega$ with coefficients in $\cO_{K_v}$.
 The method
of Coleman-Chabauty allows one to bound $|X(K)|$
 in terms of the number of zeros of the $p$-adic analytic
function $\lambda_w$.
%, a $p$-adic analytic
%function which has power series representations convergent on residue
%classes of $X(K_v)$.  

 Coleman considers the case where the  curve $X$ has good reduction at $v$, 
(that is, where
$X/K_v$ has a smooth model $\cX$
over $\cO_{K_v}$).  
In this case, some
multiple of the differential $\omega$ reduces to a differential $\wto$ on the special
fiber $\cXb$ and the zeros of $\wto$ on $\cXb$ give rise to
information about the coefficients in the power series expansion of
$\lambda_\omega$.  Coleman uses this information, along with some
simple Newton polygon arguments to obtain a variety of results on the
size of $X(K)$ in \cite{Co1} and \cite{Co2}.  For example, he shows
in \cite[0.ii]{Co2}
that, 
\begin{emp} \label{Co2-eq} If $X$ is a curve of genus $g$ defined over a number field $K$ with completion $K_v$ unramified over $\Qp$, and if  $p > 2g$ and $X$
has good reduction over $\cO_{K_v}$, then 
 $$ |X(K)| \leq q - 1 + 2g (\sqrt{q} + 1),$$
whenever $\Chab(J,K,v) < g$. 
\end{emp}

Let $\cX/\cO_{K_v}$ be any regular model 
  of $X/K_v$.
Let $\cXb_{ns}(\Fq)$ denote the 
nonsingular locus of $\cXb(\Fq)$.
In this section, we extend Coleman's result and prove:
\begin{thm}\label{to-use}
Let $X/K$ be a curve of genus $g$ defined over a number field $K$ with
completion $K_v$ unramified over $\Qp$.
  If $\Chab(J, K, v) < g$ and $p^2 > 2g+1$, then
  for any subset $\cU \subset \cXb_{ns}(\Fq)$ of the special fiber $\cXb$ of a regular model $\cX/\cO_{K_v}$
  of $X/K_v$, we have
 $$ | r^{-1}(\cU) \cap X(K) | \leq |\cU| + \left(\frac{p-1}{p-2} \right)(2g - 2).$$

\end{thm}

Thus, this theorem removes from the original  method of Chabauty-Coleman the hypothesis that $X/K$ has good reduction
at the place $v$.
McCallum pointed out at the Arizona Winter School
meeting of 1999 that
such  a generalization of the method 
 should be possible.  In fact, the early sections of \cite{Mcc1} 
can probably be
used to obtain results similar to \ref{Co2-eq} in the case where the 
curve has  a regular model with a reduced special fiber.  Before presenting the proof of
Theorem \ref{to-use}, let us
 prove a  simple lemma that allow us to bound the
number of zeros of  $\lambda_\omega$ in terms of information about
local power series expansions.  Similar arguments can be found in
\cite{Co1}, \cite{Co2}, \cite{Mcc1}, and \cite{Weth}. For simplicity,
let us assume that $K_v/{\mathbb Q}_p$ is unramified. 
Let $\lambda$ be a $p$-adic analytic function
$$\lambda: X(K_v) \lra K_v.$$
Let $P \in X(K_v)$ with reduction $r(P)
= Q$ in $\cX$ and let $u: D_Q \lra p \cO_{K_v}$ be a local coordinate at
$P$.  Suppose that $\lambda$ has a power series expansion of the form
  $$
  \lambda = a_0 + \sum\limits_{m=1}^{\infty} \frac{a_m}{m}
  u^m, $$
  where $a_m \in \cO_{K_v}$, and $v(a_m) = 0$ for some $m$. We can thus consider
$\lambda $ as a power series $\lambda(u)$ in the variable $u$, converging
on the disk $|u| \leq |p|$.
  The $p$-adic
  Weierstrass preparation theorem (\cite[Thm. 14]{Koblitz}) allows us
  to bound the number of zeros of $\lambda$ in $D_Q$.  As
  this result is most easily stated on the  disc $\cO_{K_v}$, we
  will  make the substitution $z := u/p$.  This
  gives us a power series expansion for $\lambda$ in $z$ as
  $$
  \lambda(z) = a_0 + \sum\limits_{m=1}^{\infty} \frac{a_m}{m} p^m
  z^m,$$
converging for all $z \in \cO_{K_v}$.
  Let us define
  $$I(\lambda, D_Q) := \min \{ m \mid v(a_m) = 0 
%\quad \textnormal{ and $ m \geq 1$} 
\}$$
  and
  $$
  J(\lambda, D_Q) := \min\limits \{ m \mid v(a_\ell p^\ell/\ell) >
  v(a_m p^m/m) \quad \textnormal{ for all $ \ell > m \} $}. $$
The Weierstrass preparation theorem then implies that the number of $z \in
  \cO_{K_v}$ for which $\lambda(z) = 0$ is at most $J(\lambda,D_Q)$. 
It also follows from this theorem that when $I(\lambda, D_Q)> 0$,
the number of $z \in
  \cO_{K_v}$ for which $\lambda'(z) = 0$ is at most $I(\lambda,D_Q)-1$ 
(where $\lambda'(z)$ denotes the formal derivative of $\lambda(z)$).
 It will  be convenient in the   lemma below to use the function 
  $$ \rho(x) := x - \log_p x.$$
It is clear that $\rho(m)$ is a lower bound for
$v( a_m p^m/m)$ when $a_m \in \cO_{K_v}$, since
$$ v(a_m p^m/m) = m + v(a_m) - v(m) \geq m -\log_p m.$$
The derivative of $\rho(x)$ is 
$ \rho'(x) = 1 - 1/x \ln p$. Thus, when $p > 2$,
  $\rho'(x) > 0$ for $x \geq 1$, and
 $\rho$ is therefore an increasing function for $x \geq 1$.

\begin{lem}\label{p-does} Let  $p > 2$. Assume that $K_v/{\mathbb Q}_p$
is unramified and that $I(\lambda, D_Q) < p^2  - 2$.
\begin{enumerate}[\rm a)] 
\item
 Suppose that $p \mid ( I(\lambda, D_Q) + 1 )$.  
Then $J(\lambda, D_Q) \leq I(\lambda, D_Q) + 1$.
\item Suppose that $p \nmid (I(\lambda, D_Q)+ 1)$.  Then $J(\lambda, D_Q) \leq I(\lambda, D_Q)$.
\end{enumerate} 
\end{lem}
\noindent{\bf Proof:}
To simplify notation, we will denote $I(\lambda, D_Q)$ as $I$.  Note
that $$ \rho(I+2) = I + 2 - \log_p(I+2) > I,$$
since $I + 2 < p^2$.
Let us now prove a). Then  $I>0$  
and  $v(I)=0$.
Since
$v(a_I) = 0$, we see that $v(\frac{a_I}{I} p^I) = I$.  Hence,
$$ \rho(I + 2) > v(\frac{a_I}{I} p^I).$$
Since $\rho(x)$ is increasing for $x>1$, it follows that for all $j
\geq I+2$, we have $v(a_jp^j/j) \geq \rho(j) > I$.  If
$v(\frac{a_{I+1}}{I+1} p^{I+1}) > I$, then
 $J(\lambda, D_Q) =  I(\lambda, D_Q)$, and if 
$v(\frac{a_{I+1}}{I+1} p^{I+1}) = I$,
then $J(\lambda, D_Q) =  I(\lambda, D_Q) + 1$. 

Part b) is clear when $I=0$. To prove b) when $I>0$, it is easy to see
that we need only show that 
$v(\frac{a_j}{j} p^j) > I$
for all $j > I$, since $v(\frac{a_Ip^I}{I}) \leq I$ for $I > 0$.  
Now, since $p \nmid (I(\lambda, D_Q)+ 1)$, we find  that 
$$ v(\frac{a_{I+1}}{I+1} p^{I+1}) \geq I + 1 > I.$$
Recall  that $ \rho(I+2)   > I$.
 Using the fact that $\rho(x)$ is increasing for
$x \geq 1$, we see that $\rho(j) > I$ for all $j>I$, and we are done. 
 
\begin{remark}
Coleman uses arguments similar to the lemma above combined with
information about the terms $I(\lambda_\omega, D_Q)$ to obtain
\ref{Co2-eq}.  
In \cite{Co1} and
\cite{Co2}, he controls the terms $I(\lambda_\omega,
D_Q)$ by counting the zeros of the pull-back $\wto$ of $\omega \in
\Gamma(\cX, \Omega_{\cX/\cO_{K_v}})$ to the special fiber
$\overline{\mathcal X}$.  Indeed, when the curve has good reduction, the
differential $\wto$ must have exactly $2g(X) -2$ zeros, counted with
multiplicity, since it is a nonzero differential on a smooth
irreducible curve.  On a general regular
model, however,   $\wto$ may vanish on entire components
of the special fiber, which makes it difficult to count the zeros of
$\wto$ in a sensible manner.  Furthermore, when $\cXb$ is not reduced,
its dualizing sheaf may not even be a line bundle.  Hence, we will not
work with the reduction of $\omega$, but will proceed as follows.
%rather with certain multiples of $\omega$.
\end{remark}

Let $\cX/\cO_{K_v}$ be any regular model of $X/K_v$. 
 Let $Q \in \cXb_{ns}(\Fq)$.
Denote by  $\cO_Q$ the local ring $
\cO_{\cX,Q}$. Let $P_0 \in X(K_v)$ be a point reducing to $Q$.
The closure of $P_0$ in $\cX$ corresponds to a prime ideal of height
$1$ in $\cO_Q$. Since $\cO_Q$ is regular and, thus, a UFD, we can write
this prime ideal as $(u)$ for some prime $u \in \cO_Q$. It follows that the maximal ideal of $\cO_Q$ is generated by $\pi $ and $u$. Let
$\hat{\cO}_Q$ denote the completion of the ring $\cO_Q$ at the prime $(u)$.
One easily shows that the natural map from
the ring $\cO_{K_v}[[u]]$ of formal power series to the ring $\hat{\cO}_Q$
(which sends $u$ to $u$) is an isomorphism.
It is also easy to check that the $\hat{\cO}_Q$-module of relative differentials
$\Omega_{\hat{\cO}_Q/\cO_{K_v}}$ is generated by $du$. Any differential
$w \in \Omega_{\cO_Q/\cO_{K_v}}$ can thus be written as a power
series $\omega = \sum_{m=0}^\infty a_{m+1}  u^m du$
with $a_m \in \cO_{K_v}$ for all $m \in \bZ_{\geq 0}$.

 Let us  fix a generator $\pi$ for the
maximal ideal of $\cO_{K_v}$.
Since $Q$ is a
nonsingular point of $\cXb$, the local ring $\cO_{\cXb, Q}$ is a
discrete valuation ring, and we denote by  $v_Q$ its valuation.
For any $P \in X_{K_v}$, we denote by $v_P$ the valuation of the local ring
$\cO_{X_{K_v},P}$.  
 A differential $\omega \in \Gamma(\cX, \Omega_{\cX/\cO_{K_v}})$
 pulls back to a differential $i^*
\omega$ via the natural map $i:\cX_{K_v} \lra \cX$ from the generic
fiber $\cX_{K_v}$ of $\cX$  to $\cX$. We denote by $(i^* \omega)_0$
the divisor of zeros of $i^* \omega$, and we shall write $(i^* \omega)_0=
\sum_P v_P(i^* \omega) P$.
 We have the following
proposition relating the power series expansion of some multiple of
$\omega$ to the zeros of $i^* \omega$ in $r^{-1}(Q)$.  

\begin{prop}\label{CC-reg} Keep the notation introduced above.
  Let $\omega \in \Gamma(\cX, \Omega_{\cX/\cO_{K_v}})$ and let $Q \in \cXb_{ns}(\Fq)$.
  Then there exists an element $t \in  K_v$
such that $t \omega \in \Gamma(\cX, \Omega_{\cX/\cO_{K_v}})$ 
and has a local power
series expansion
\begin{equation}\label{t-omega1}
 t \omega = \sum\limits_{m=0}^\infty a_{m+1}  u^m du,
\end{equation}
with $a_m \in \cO_{K_v}$ for all $m \in \bZ_{\geq 0}$, such that
\begin{equation}\label{t-omega2}
\min \{ m \mid v(a_{m+1}) = 0 \} = \sum_{  r(P) = Q  }  [K_v(P) : K_v] v_P(i^* \omega),
\end{equation}
where the sum is taken over all points $P$ of the scheme $X_{K_v}$
such that the intersection of the closure of $P$ in $\cX$ with $\cXb$
is $Q$.
\end{prop}

\noindent{\bf Proof:} Since $\Omega_{\cX/\cO_{K_v}}$ is locally free 
of rank $1$, let $f$ denote a generator
of the stalk $\Omega_{\cO_Q/\cO_{K_v}}$ of 
$\Omega_{\cX/\cO_{K_v}}$ at $Q$. We
can write the stalk of $\omega$ at $Q$ as $s f$, where $s \in \cO_Q$.   Factor $s$ as
$$
s = \gamma_1^{\ell_1} \cdots \gamma_n^{\ell_m} \pi^{\ell'}$$
where
the $\gamma_j$ are generators for primes corresponding to points $P_j$
on the generic fiber of $\cX$.  It is not hard to see that
$v_{P_j}(i^* \omega) = \ell_j$. Indeed, one obtains the local ring
$\cO_{X_{K_v},P_j}$ localizing $\cO_Q$ at the prime ideal generated by
$\gamma_j$, so we see that the ideal generated by $s$ in $\cO_{X_{K_v},P_j}$ 
is just $\cm_{P_j}^{\ell_j}$, where $\cm_{P_j}$ is the maximal
ideal in $\cO_{X_{K_v},P_j}$; since $f$ pulls back to a generator for the
stalk of $\Omega_{{\cX_{K_v}}/K_v}$,  $s f$ must pull back to a
differential with order of vanishing equal to $v_{{P_j}}(s) =
\ell_j$ for all $j$.

After dividing $s$ by $\pi^{\ell'}$ we obtain an element $s_1$ that is
not in $\pi\cO_Q $.  Now complete  $\cO_Q$  at $(u) $. We obtain a
power series expansion for $s_1f$:
$$ s_1f = \sum\limits_{m=0}^\infty a_{m+1} u^m du.$$
Let us denote by $s_2$ the element $\sum_{m=0}^\infty a_{m+1} u^m$ of $\hat{\cO}_Q$. The elements $s_1$ and $s_2$ differ by a unit in $\hat{\cO}_Q$.
Consider the commutative diagram 
$$
\begin{CD}
\cO_Q @>>> \hat{\cO}_Q\\
@VVV @VVV\\
\cO_Q/(\pi) @>>> \hat{\cO}_Q/(\pi).
\end{CD}
$$
It is easy to check that $\hat{\cO}_Q/(\pi)$ is the completion
of $\cO_Q/(\pi)$ at the maximal ideal $(u)$. Thus the valuation
$v_Q$ of $\cO_Q/(\pi)$ extends to a valuation  on $\hat{\cO}_Q/(\pi)$, again denoted by $v_Q$.
Denoting by $\phi_\pi$ the map taking $\hat{\cO}_Q$ to
$\hat{\cO}_Q/(\pi) $, it is clear that
$$
\min \{ m \mid v(a_m) = 0 \} = v_Q(\phi_\pi(s_2)). $$
Since
$$ v_Q(\phi_\pi(s_2)) =  \sum\limits_{j=1}^m \ell_j v_Q(\phi_\pi(\gamma_j)),$$
it will suffice to show that
$$  v_Q(\phi_\pi(\gamma_j)) = [K_v(P_j):K_v]. $$
This follows from the fact that 
\begin{equation}
\begin{split}
v_Q(\phi_\pi(\gamma_j)) & = \dim_{\Fq} \big( (\cO_Q/\pi \cO_Q)/\phi_\pi(\gamma_j)
\big) \\
& = \rank_{\cO_{K_v}}  (\cO_Q/\gamma_j \cO_Q) \\
& =  [K_v(P_j):K_v],
\end{split}
\end{equation}
since $(\cO_Q/\gamma_j \cO_Q)$ is a free $\cO_{K_v}$-module (which follows from
the fact that $\cO_{K_v}$ is of course a principal ideal domain). This concludes the proof of  \ref{CC-reg}.

\medskip
Let us now apply Lemma \ref{p-does}  and
Proposition \ref{CC-reg} to the sort
of $p$-adic analytic function that arises in the Chabauty-Coleman
method. 

\begin{prop}\label{prop-to-use}
Let $X/K$ be a curve of genus $g$ defined over a number field $K$ with completion $K_v$ unramified over $\Qp$.
 Let $\cX/\cO_{K_v}$ be any regular
model for $X/K_v$, and let $\cU \subset \cXb_{ns}(\Fq)$.  If
$p^2 > 2g+1$ and $\lambda_\omega$ is as in \eqref{eq.lambda}, then
$$| r^{-1}(\cU) \cap \lambda_\omega^{-1}(0) | \leq |\cU| + \left( 
\frac{p-1}{p-2} \right)(2g - 2).$$
\end{prop}

\noindent{\bf Proof:}
  Choose $Q \in \cU$.  For any nonzero element $t \in
  K_v$, multiplying $\lambda_\omega$ by $t$ will not change the zeros
  of $\lambda_\omega$, and, furthermore, $\lambda_{t \omega} = t
  \lambda_\omega$ , so
$$ | r^{-1}(Q) \cap \lambda_\omega^{-1}(0) | = | r^{-1}(Q) \cap
\lambda_{ t \omega}^{-1}(0) | . $$
%%Tom3 removed  \lambda_{ t \omega,D}^{-1}(0)
  Thus, we may choose $t \in K_v$ and apply Proposition
  \ref{CC-reg} to obtain a power series expansion of the form
  \eqref{t-omega1} for which equation \eqref{t-omega2} holds.
We denote by $Z(\omega, Q)$ the sum
$$\sum_{ r(P)=Q }
%\begin{Sb} r(P) = Q \\ 
%P \in \cX({\overline K}_v) \end{Sb}}  
[K_v(P) : K_v] v_P(i^* \omega)$$
appearing in the statement of Proposition \ref{CC-reg}.
It is clear that since $d(\lambda_{tw}) = tw$,   $$Z(\omega, Q) = I(\lambda_{t \omega},D_Q) - 1.$$ 
Since
$$\sum\limits_{Q \in \cU} Z(\omega, Q) \leq
\sum\limits_{ P \in X_{K_v}} [K_v(P) : K_v]v_P(i^*\omega)  = 2g - 2 < p^2-3,$$
we find that $I(\lambda_{t \omega},D_Q) < p^2-2$, and we can apply 
Lemma \ref{p-does}. We obtain
\begin{equation}\label{vQ-bound}
\begin{split}
  | r^{-1}( \cU ) \cap \lambda_\omega^{-1}(0) | & \leq 
\sum_{Q \in \cU} J(\lambda_{tw}, D_Q)\\
& \leq  \sum\limits_{p \mid  
(Z(\omega, Q)+2)} (Z(\omega, Q) + 2) + \sum\limits_{p \nmid
  (Z(\omega, Q) + 2) } (Z(\omega, Q) + 1).
\end{split}
\end{equation}
 If $p \mid (Z(\omega, Q) + 2)$, then
 $Z(\omega, Q) \geq p - 2$.
Since
$\sum_{Q \in \cU} Z(\omega, Q) \leq 2g - 2$,
  there are at most $(2g - 2)/(p-2)$
points $Q \in \overline{\cX}_{ns}(\overline{\mathbb F}_q)$ 
 for
which $p \mid (Z(\omega, Q) + 2)$.  Hence, \eqref{vQ-bound} becomes
\begin{equation*}
\begin{split}
  | r^{-1}( \cU ) \cap \lambda_\omega^{-1}(0) | & \leq \sum\limits_{Q
    \in \cU}
  Z(\omega, Q) + | \cU | + \frac{2g-2}{p-2} \\
  & \leq | \cU | + \left(1 + \frac{1}{p-2} \right)(2g - 2).
\end{split}
\end{equation*}
We are now ready prove Theorem \ref{to-use}.

\begin{emp} {\bf Proof of \ref{to-use}:}
  Each differential $\eta \in \Gamma(J, \Omega_{J/K_v})$ gives rise
  to a homomorphism
  $ \lambda_{\eta}: J(K_v) \lra K_v$. 
  Since
 $\Chab(J,K,v)  < {\rm dim}H^0(J, \Omega_{J/K_v})$,
  there must be a nonzero $\eta$ for which
  $ \lambda_\eta (J(K)) = 0$.
   We may assume that
  $X(K)$ contains a point $Q$, as otherwise our assertion is
  trivial.  Hence, we may embed $X(K_v)$ into $J$ via the mapping
  $ j: X \to J$, which sends $P \in X(K_v)$ to the class of $ P-Q$.
   Now, $\eta$ pulls back to a
  differential $\omega$ on $X$, and $\lambda_\eta$ restricts to a
  function $\lambda_\omega$ that vanishes on $X(K)$ (because $j$ sends
  points in $X(K)$ to points in $J(K)$).  Applying Proposition
  \ref{prop-to-use} then gives the desired result.
\end{emp}

\begin{emp} \label{tobedone}
Note that if an abelian variety $A/K$ is $K$-isogenous to a product
$\prod A_i$, then $\Chab(A,K,v) = \sum \Chab(A_i,K,v)$. Thus, the method
of Chabauty-Coleman can be applied to $A$ if and only if $\Chab(A_i,K,v)
< \dim(A_i)$ for some $i$.
\end{emp}

When the Chabauty rank of $X/K$ is zero, we can strengthen 
Theorem \ref{to-use} as follows.

\begin{prop} \label{pro.00}
  Let $X/K$ be a curve of genus $g$ defined over a number field $K$
  with completion $K_v/ {\mathbb Q}_p$ unramified over $\Qp$.
Let $\cX/\cO_{K_v}$ be any regular model for $X/K_v$. 
   If $\Chab(J,K,v) = 0$,
then $ |X(K)| \leq |\cXb_{ns}(\Fq)| $.
\end{prop}
\noindent{\bf Proof:}
We claim that for each $Q \in \cXb_{ns}({\mathbb F}_p)$, 
the set $r^{-1}(Q)$
contains at most one $K$-rational point of $X$.  We begin by noting
that since $\Chab(J,K,v) = 0$, every differential $\eta \in \Gamma(J,
\Omega_{J/K_v})$ has the property that its formal $v$-adic integral
$\lambda_\eta$ is identically 0 on $J(K)$.  It follows that for every
differential $\omega \in \Gamma(X, \Omega_{X/K_v})$, the
formal integral $\lambda_\omega$ is identically 0 on $X(K)$.  
To prove our claim, we need only to show that
for some
$w \in \Omega_{\cX/\cO_{K_v}}$, the sum
$$
\sum_{ r(P) = Q } [K_v(P) : K_v] v_P(i^* \omega) $$
appearing in the statement of \ref{CC-reg} is equal to $0$.
Indeed, the  proof of Proposition \ref{prop-to-use} shows that 
this sum is equal to $I(\lambda_{w}, D_Q) -1$, and that
$$ |D_Q \cap \lambda_w^{-1}(0) | \leq J(\lambda_w, D_Q) \leq I(\lambda_{w}, D_Q).$$
  Now, since
the sheaf $\Omega_{\cX/\cO_{K_v}}$ is generated by global sections,
there must be a global section $\omega \in \Gamma(\cX,
\Omega_{\cX/\cO_{K_v}})$ whose stalk at $Q$ generates the stalk
$\Omega_{\cX/\cO_{K_v},Q}$ of $\Omega_{\cX/\cO_{K_v}}$ at $Q$ as an
$\cO_{Q,\cX}$-module.  Otherwise, the stalks of all the global
sections would be in $\mathcal{M}_Q \Omega_{\cX/\cO_{K_v}, Q}$ and
hence could not together generate $\Omega_{\cX/\cO_{K_v}, Q}$ as an
$\cO_{K_v}$-module.  
%Applying Proposition \ref{CC-reg} and 
%observing
Since $\omega$ has to generate the stalks of the differential sheaf at
all of the $P$ with $r(P) = Q$,
we find that the sum
$
\sum_{ r(P) = Q } [K_v(P) : K_v] v_P(i^* \omega) $ appearing in \ref{CC-reg}
must be 0. 
This concludes the proof of \ref{pro.00}.

\begin{remark}
When the Mordell-Weil rank of $X/K$ is zero, Proposition \ref{pro.00}
can be strengthened as follows.
 Let $X/K$ be a curve of genus $g$ defined over a number field $K$ with completion $K_v/ {\mathbb Q}_p$ such that $v(p) < p-1$.
 Let $\cX/\cO_{K_v}$ be any regular
model for $X/K_v$.   If the Mordell-Weil rank of $X/K$ is zero, then
$|X(K)| \leq |\cXb_{ns}(\Fq)| $.

To prove this statement, we proceed as follows.
We claim that for each $Q \in \cXb_{ns}(\Fq)$, the set $r^{-1}(Q)$
contains at most one $K$-rational point of $X$.
Indeed, suppose that $P$ and $P'$ belong to $r^{-1}(Q) \cap X(K)$.
Consider the embedding of $X$ in $J$ using the map $D \mapsto D- \deg(D)P$.
The map $X \to J$ extends to a map from the smooth part of $\cX $ 
to the N\'eron model
of $J$.  By construction,  $P'-P$ reduces to the origin. In other words,
$P'-P$ belongs to the kernel of the reduction, which does not contain any 
torsion point (see for instance \cite{Ser}, LG 4.25-4.26). Thus, $P'-P$ has infinite order, contradicting the assumption that
the Mordell-Weil rank over $K$ is zero. 
\end{remark}

The following generalization of  
\ref{Co2-eq} will not be used in this paper. 

\begin{cor} \label{cor.g-1}
Let $X/K$ be a curve of genus $g$ defined over a number field $K$ with completion $K_v$ unramified over $\Qp$.
 Let $\cX/\cO_{K_v}$ be any regular
model for $X/K_v$.  If
$p > 2g$ and $\Chab(J,K,v) < g$, then
$|X(K)| \leq |\cXb_{ns}(\Fq)| + 2g-2 $.
\end{cor}

\noindent{\bf Proof:}
  We proceed exactly as in \ref{to-use}, with $\cU =\cXb_{ns}(\Fq)$. 
   Since $p>2g$, we find  that $(2g-2)/(p-2) < 1$, and $|X(K)|$ is of
  course a whole number. So
  $|X(K)| \leq |\cXb_{ns}(\Fq)| + (2g - 2)$,
  as desired.
 
\medskip In view of \ref{pro.00} and \ref{cor.g-1}, it is natural to wonder,
under the hypotheses of \ref{cor.g-1}, whether the bound for $|X(K)|$
can be made to depend on the precise value of $\Chab(J,K,v)$, such as
a bound of the form $|X(K)| \leq |\cXb_{ns}(\Fq)| +2 \Chab(J,K,v)$.

\section{Constructing regular models of curves}

Let $K$ be a field with a discrete valuation $v_K$. Let ${\mathcal
  O}_K$ denote the ring of integers of $K$, with maximal ideal
$(\pi_K)$ and residue field $k$. Let $p:= {\rm char}(k)$.  When no
confusion may ensue, we will write $\pi$ and $v$ instead of $\pi_K$
and $v_K$.  Let $X/K$ be the nonsingular proper model of the plane
curve $C/K$ given by a homogeneous equation $f(x, y, z) \in {\mathcal
  O}_K [x, y, z]$.  Assume that the ideal of ${\mathcal O}_K$
generated by the coefficients of $f$ is ${\mathcal O}_K$.  A regular
model ${\mathcal X}/{\mathcal O}_K$ of $X/K$ can be theoretically
obtained as follows.  Consider first the model ${\mathcal D}/{\mathcal
  O}_K$ given by
$$
{\mathcal D} := {\rm Proj}({\mathcal O}_K[x, y, z]/(f)).
$$
and let ${\mathcal X}$ be the minimal desingularization of the
normalization ${\mathcal D}^{nor}$ of ${\mathcal D}$.  Let $\rho:
{\mathcal X} \rightarrow {\mathcal D}^{nor}$ denote the
desingularization map and let $\nu: {\mathcal D}^{nor} \rightarrow
{\mathcal D}$ be the normalization map.  Both $\rho$ and $\nu$ are
usually quite difficult to describe explicitly, even in the case of a
rather simple equation $f$.  Obviously, we have the option of changing
the defining equation $f$, but this does not always lead to more
easily described maps $\rho$ and $\nu$.

Another way of constructing a regular model ${\mathcal X}/{\mathcal O}_K$ of $X/K$
was introduced by Viehweg in \cite{Vie}.  
It may happen that over a Galois extension $L/K$, a normal model
${\mathcal Y}/{\mathcal O}_L$ of $X_L/L$ can be described.  If the Galois group ${\rm
Gal}(L/K)$ acts on ${\mathcal Y}$, lifting its action on ${\rm Spec}({\mathcal O}_L)$,
then we may consider the quotient ${\mathcal Y}/{\rm Gal}(L/K)$ as a scheme over
${\rm Spec}({\mathcal O}_K)$.  The scheme ${\mathcal Y}/{\rm Gal}(L/K)$ is a normal
model of $X/K$ and, thus, a desingularization $\rho: {\mathcal X} \rightarrow {\mathcal
Y}/{\rm Gal}(L/K)$ leads to a regular model ${\mathcal X}/{\mathcal O}_K$ of $X/K$.
A key feature of this method is the fact that
the singularities of ${\mathcal Y}/{\rm Gal}(L/K)$ are quotient singularities
and that when $L/K$ is tame, such singularities are well-understood.

In this section, we first use this method to construct regular models
of curves having potentially good reduction after a tame extension
$L/K$, such as the superelliptic curves $X:= X_{F,h}$ with $\pi \nmid
d^*(F)$ and $p>n$ (see \ref{lem.tame} below).  The model ${\mathcal
  Y}/{\mathcal O}_L$ that we use in this case is simply the smooth
minimal model of $X_{L}/L$, where $L/K$ is large enough to ensure that
$X_{L}/L$ has good reduction.
 
We shall also use the above method to consider the more difficult
case where
$\pi \mid d^*(F)$ and $\pi \mid h$.
In this case, we construct just enough  of a
regular model for $X/K$ to be able to bound the number of residue
classes of primitive integral solutions to the Thue equation $F(x,y) =
h$.  Here the field $L/K$ will be the 
splitting field over $K$ of the polynomial $F(x,1)$,
and the model   ${\mathcal Y}/{\mathcal O}_L$ will be  the normalization
of the model 
$${\mathcal C} := \Proj {\mathcal O}_L[x,y,z]/(hz^n - F(x,y)).$$
Some smooth open affine subsets of the model ${\mathcal Y}$ are described 
below in \ref{emp.smoothopen}.
Let $F(x, 1) = c \prod^s_{i=1}(x - \alpha_i)^{n_i}$ in 
$\overline{K}[x]$.  Let 
$$d^*(F) :=
c \prod_{i\ne j}(\alpha_i - \alpha_j) \in \overline{K}.$$
We shall say that $\pi_K \nmid d^*(F)$ if $\pi_K \notin (d^*(F))$ in 
${\mathcal O}_{\overline{K}}$.

\begin{lemma}  \label{lem.tame} Assume that ${\rm char}(k) \nmid n$.
Let $X:= X_{F,h}/K$.
\begin{enumerate}[\rm a)]
\item
If $\pi_K \nmid d^*(F)$ and $\pi_K \nmid h$, then $X/K$ has good
reduction.
\item If $\pi_K \nmid d^*(F)$ and $\pi_K \mid h$, then $X/K$ achieves good
reduction over $L:= K(\sqrt[n]{h})$.
%, and 
%$X_M/M$ does not have good reduction if $  M$ does not contain $L$.
\end{enumerate}
\end{lemma}

\noindent{\bf Proof:} Consider the model
${\mathcal C}/{\mathcal O}_K$ given by ${\rm Proj}({\mathcal O}_K[x,y,z]/(hz^n - F(x, y))$ and its
normalization ${\mathcal C}^{nor}/{\mathcal O}_K$.  The generic fiber of ${\mathcal
C}^{nor}$ has genus equal to
$$
2g(X) - 2 = n(s -2) - \sum^s_{i=1}{\rm gcd}(n, n_i).
$$
If $g(X) = 0$, then $X/K$ has obviously good reduction over 
${\mathcal O}_K$.
Assume then that $g(X) > 0$.
Since $\pi_K \nmid h$,
the reduction modulo $\pi$ of $hz^n - F(x, y)$ is irreducible. When
$\pi_K \nmid h$ and $\pi_K \nmid d^*(F)$, we find that the geometric genus of
${\mathcal
C}^{nor}_k$ is equal to $g(X)$.  Thus, ${\mathcal
C}^{nor}_k$ is non-singular since its arithmetical genus
is equal to the genus of $X$. It follows that
 ${\mathcal C}^{nor}/{\mathcal O}_K$ is the
(minimal) regular model of $X/K$.

If $\pi_K \mid h$, consider the change of variables $z'= \sqrt[n]{h} z,
x'=x $ and $y'=y$. Then ${\rm Proj}({\mathcal O}_L[x',y',z']/({z'}^n-F(x',y')))$
is a model for $X_{L}/L$. Hence, we may apply a) to find
that $X_{L}/L$ has good reduction. 
%The proof of the second statement
%of b) will be given in \ref{rem.lateproof}. 
This concludes the proof of \ref{lem.tame}.
\medskip

For most of the applications that we have in mind,  the residue
field ${\mathcal O}_K/(\pi_K)$ will be ${\mathbb F}_p$, and 
we will assume that $p>n$. The following lemma shows that
we may assume, under these hypotheses, that $F(x,1)$ is monic.

\begin{lemma} Assume that $|{\mathcal O}_K/(\pi_K)| > n$.
Then, up to a change of variable, we may assume that $F(x,1)$ is monic in ${\mathcal O}_K[x]$.  
\end{lemma}
\noindent{\bf Proof:} Write $F(x,y) = \prod\limits_{i=1}^{n}(\beta_i x - \rho_i y)$, with
$\beta_i, \rho_i$ in ${\mathcal O}_K$. By substituting $y' = y + ux$,
we can write
$$F(x,y') = \prod\limits_{i=1}^n ((\beta_i + \rho_i u)x - \rho_i y).$$
Since the
coefficients of $F$ have no common factor, we must have
$\min(v(\beta_i), v(\rho_i)) = 0$ for each $i$, so we can
choose $u \in {\mathcal O}_K$ such that $ \rho_i u + \beta_i   \in {\mathcal O}_K^*$
for all $i$ by simply avoiding the residue classes of the
$-(\beta_i/\rho_i)$ for which $v(\rho_i ) =0$.  There are at
most $n$ residue classes to be avoided and there are at least $n + 1$
residue classes in ${\mathcal O}_K/(\pi_K)$ by hypothesis, 
so we can do this.  This
allows us to rewrite our original equation $F(x,y) = h$ as
$$\prod_{i=1}^{n}(x - \alpha_i y) = \mu \pi^w,$$
where $\mu$ is unit in ${\mathcal O}_K$ and $w = v(h)$.
\medskip

\begin{emp} {\bf Some regular affine subsets of the normalization of ${\mathcal C}$.} \label{emp.smoothopen}
\end{emp}

In what follows, we assume that $F(x,1)$ is monic, 
 that $ \pi_K \mid  h$ and that $\pi_K \mid d^*(F)$.
Assume also that $K$ is complete, so that for any finite extension $L/K$,
the integral closure ${\mathcal O}_L$ of ${\mathcal O}_K$ in $L$
is a local ring. 

Let $L/K$  be the 
splitting field over $K$ of the polynomial $F(x,1)$.
We denote by $v$ the valuation of ${\mathcal O}_L$,
and let $\pi$ be a uniformizer of  ${\mathcal O}_L$.
Let us say that $P=(a,b) $ is a {\it primitive integral solution} to $F(x,y)=h$
if $a,b \in {\mathcal O}_K$ and $\gcd(a,b)=1$.
We describe below an affine regular scheme ${\mathcal U}/{\mathcal
  O}_L$ such that ${\mathcal U} \times_{{\rm Spec}({\mathcal O}_L)}
{\rm Spec}(L)$ is open in $X_{F,h,L}$ and $P \in {\mathcal U}(L)$ has
a non-trivial reduction modulo $(\pi)$. In other words, the closure of
$P$ in ${\mathcal U}$ includes a point on the special fiber  of $\cU$.

Consider any root $\alpha_i$ of $F(x,1)$ such that  $v(a -
\alpha_i b) = \max_j(v(a - \alpha_j b))$. Let
$t:= v(a -
\alpha_i b)$.
  Change variables from $x$ to $z_0:=
x- \alpha_iy $, so that $F(z_0,y) = z_0\prod_{j=1}^{n-1}(z_0
-\gamma_j y)$, where $\gamma_j := \alpha_j -\alpha_i$ for $j < i$ and
$\gamma_j :=\alpha_{j+1} - \alpha_i$ for $j \geq i$. Define $s_0 := 0$, and
then recursively define
  $$ s_k := \min \{ v(\gamma_j) \mid t \geq v(\gamma_j) > s_{k-1}
  \},$$
  for $k \geq 1$.  We obtain in this way a finite increasing
sequence of integers.
If $t$ is not the largest integer of this sequence, add $t$
to the sequence. Denote the elements of the new sequence by
  $$s_0 < s_1 < \cdots < s_m = t.$$
 Define, for $k<m$,
  $$\cS_k := \{ \gamma_j \mid v(\gamma_j) = s_k \}.$$
The set $S_m$ is defined to be $\{ \gamma_j \mid v(\gamma_j)  \geq s_m \}$.
If $\gamma$ is a root of $F(z_0,y)$,
let $n(\gamma)$ denote its multiplicity.
Then, for $k \leq m$,  define
  $z_k$ to be $z_0/\pi^{s_k}$, and   let
     $F_k$ be the polynomial
  $$F_k(z_k, y) := 
\prod_{j=0}^{k} \prod_{\gamma \in {\mathcal S}_j}(\pi^{s_k-s_j} z_k 
- \gamma \pi^{-s_j}y)^{n(\gamma)} 
\prod_{\gamma \notin \cup_{j=0}^k {\mathcal S}_j}
(  z_k 
- \gamma \pi^{-s_k}y)^{n(\gamma)}.$$
Set
$$u_k := \sum_{j=0}^k  (\sum_{\gamma \in {\mathcal S}_j} n(\gamma)  )s_j + \sum_{j=k+1}^m  (\sum_{\gamma \in {\mathcal S}_j} n(\gamma) )s_k.$$ 
Then $F_k(z_k,y) = F_0(z_0,y)  \pi^{-u_k}$. 
Finally,   let
$$ A_k := {\mathcal O}_L [z_k,y]/(F_k(z_k,y) - \mu \pi^{w -u_k}),$$
for $k \leq m$ (recall that $h = \mu \pi^w$). 
Note now that when $(a,b)$ is primitive and $\pi \mid h$, 
then $v(b) = 0$.
Indeed, if $v(b)>0$ and $(a,b)$ is primitive, then $v(a) = 0$.
Thus, $v(a-\alpha_j b) = 0$ for all $j$, contradicting
the fact that $v(F(a,b)) =v(h) >0$. Hence, $v(b)=0$.
If follows that for $j \neq i$, the inequality
$$v(a-b\alpha_j) \geq {\rm min}(v(a-b\alpha_i), v(b\alpha_i-b\alpha_j))$$
implies that either $v(a-b\alpha_j) =t$ and 
$v(\alpha_i-\alpha_j) \geq t$, or $v(a-b\alpha_j) <t$ and 
$v(a-b\alpha_j) =v(\alpha_i-\alpha_j)$.
In particular, we find that when $(a,b)$ is primitive,
$w=u_m$.

\begin{lem}
The ring $A_m$ is regular and $\Spec(A_m)/\Spec({\mathcal O}_L)$ is smooth.
\end{lem}
\noindent{\bf Proof:}
The generic
fiber of $A_m$ is easily checked to be smooth.
 Hence, we need only check points on
the special fiber of $A_m$.  We note that modulo $\pi$, the
equation $F_m(z_m,y) = \mu$ is equivalent to the equation
%\begin{equation}%\label{barthing}
$$
\left( \prod\limits_{\gamma \in {\mathcal S}_m} \left(z_m -
  (\overline{\gamma/\pi^t})y \right)^{n(\gamma)}\right) \left( \prod\limits_{j=0}^{m-1}
\prod\limits_{\gamma \in {\mathcal S}_j}
  \left(- (\overline{\gamma/\pi^{s_j}})y \right)^{n(\gamma)} \right) - \overline{\mu} = 0,
$$
%\end{equation}
because $\pi^{w -
u_m} = 1$ as noted earlier.  Since $\overline{\mu} \neq 0$, this
equation defines a nonsingular affine curve in ${\mathbb A}^2$.  Thus, the
special fiber of $A_m$ is nonsingular; therefore all the points on the
special fiber of $A_m$ are regular and $A_m$ is a regular ring.
%\end{proof}

\medskip
Let ${\mathcal Y}/{\mathcal O}_L$ denote the normalization of the scheme 
$${\mathcal C} := \Proj {\mathcal O}_L[x,y,z]/(hz^n - F(x,y)).$$

\begin{lem} \label{lem.open}
The affine scheme $\Spec A_m$ is an open subset of ${\mathcal Y}$.
\end{lem}
\noindent{\bf Proof:} There is a natural ring homomorphism
$A_{k-1}\to A_k$ that sends $z_{k-1}$ to $\pi^{s_k - s_{k-1}}z_{k}$.
Define
$$G_k(z_k,y) :=\prod_{j=0}^{k} 
\prod_{\gamma \in {\mathcal S}_j}(\pi^{s_k-s_j} z_k 
- \gamma \pi^{-s_j}y)^{n(\gamma)}.$$
Let $S_k$ denote the multiplicative subset of $A_k$
generated by $G_k(z_k,y)$.
We claim  that
$A_k$ is integral over $S_{k-1}^{-1}(A_{k-1})$.
Indeed, it suffices to show that $z_k$ is integral over
$S_{k-1}^{-1}(A_{k-1})$.
Recall that in $A_k$,
$$ F_k(z_k,y) -\mu \pi^{w-u_k}
= G_{k-1}(z_{k-1},y) \prod_{j=k}^m \prod_{\gamma \in {\mathcal S}_j}(z_k - \gamma \pi^{-s_k}y)^{n(\gamma)} -\mu \pi^{w-u_k} =0.
$$
Thus, the image of $z_k$ in $A_k$ is the
root of a monic polynomial over $S_{k-1}^{-1}(A_{k-1})$
(since $G_{k-1}$ is of course a unit in this ring).  
Hence, it follows that the map $\Spec(A_k) \to \Spec(A_{k-1})$
is quasi-finite for any $k \geq 1$. Since $A_m$ is
normal  because it is regular,
we have a natural map $j:\Spec(A_m) \to {\mathcal Y}$, and our discussion above
shows that this map is quasi-finite.
Since ${\mathcal Y}$ is normal and $j$ is generically an isomorphism, we can apply Zariski's Main Theorem to find that $j$ is an open immersion.
\if false
Recall that $A_m$ is
normal, since it is regular. We have shown that $\Spec A_m$
is integral over $\Spec(S_{m-1}^{-1}(A_{m-1}))$. Since 
$\Spec(S_{m-1}^{-1}(A_{m-1}))$ is open in $\Spec(A_{m-1})$,
we find that $\Spec A_m$ is a special open set in the normalization
of $\Spec(A_{m-1})$.  We also showed that
the normalization   of $\Spec (A_k)$ is an
open subset of the normalization   of $\Spec
(A_{k-1})$.  Since an open subset of an open
subset is itself open, we see that the $\Spec A_m$ is open subset
of the normalization of $\Spec (A_0)$. 
Furthermore, $\Spec (A_0)$ is an open
subset of $\cC$, so $\Spec A_m$ is open in ${\mathcal Y}$, as desired.  
\fi

\begin{emp} \label{emp.ualpha}
Let ${\mathcal U}(\alpha_i) := \Spec(A_m)$.
The primitive point $P=(a,b)$ in $X_{F,h}(L)$
corresponds to the point $(\pi^{-t}(a-\alpha_ib), b)$ in ${\mathcal U}(\alpha_i)(L)$.
Since this point  is integral, it has a non-trivial
reduction in the special fiber of ${\mathcal U}(\alpha_i)$.
Denote by $\cP$ the set of integral primitive solutions, so that 
$$ \cP := \{ (x,y) \in ({\mathcal O}_K)^2 \mid F(x,y) = h \mbox{ \rm   and $ \gcd(x,y) = 1 \}$}. $$
Our next lemma shows that the closure of ${\mathcal P}$ in 
the normalization ${\mathcal Y}/{\mathcal O}_L$   of  
${\mathcal C}/{\mathcal O}_L$ is contained in at most $n$ regular affine open sets;
namely, this closure is contained in the union of the 
sets ${\mathcal U}(\alpha_i)$, where $\alpha_i$ runs through
all the roots of $F(x,1)$ such that
there exists a primitive point $(a,b)$ with $v(a-\alpha_ib) = \max_j(v(a-\alpha_jb))$.
\end{emp}

\begin{lem}\label{decomp} Let $\pi \mid h$.
Let $(a,b)$ and $(a',b')$ be elements of $\cP$.   
Suppose that $v(a - \alpha_i b) = \max_j(v(a - \alpha_j
b))$ and $v(a' - \alpha_i b') = \max_j(v(a' -
\alpha_j b'))$.  Then
$v(a - \alpha_i b) = v(a' - \alpha_i b')$.
\end{lem}
\noindent{\bf Proof:}
 Recall that $v(b)=0$ when $\pi  \mid h$ and $(a,b)$ is primitive.
Suppose   that $v(a - \alpha_i b)$ and $v(a'
- \alpha_i b')$ are not equal.  We may assume without loss of
generality that
$v(a/b - \alpha_i ) > v(a'/b' - \alpha_i)$.
We claim  that this inequality implies that 
$v(a/b - \alpha_j ) \geq  v(a'/b' - \alpha_j)$ for all $j$.
Indeed, $v(a/b - \alpha_j ) \geq  v(a'/b' - \alpha_j)$ is clear
if $v(a'/b' - \alpha_i ) \leq v(a/b - \alpha_j)$. Thus we may assume
that $v(a'/b' - \alpha_i ) > v(a/b - \alpha_j)$. 
>From $v(a/b - \alpha_i ) > v(a'/b' - \alpha_i)$ we find that
$v(a/b-a'/b')=  v(a'/b' - \alpha_i)$. It follows   
from $v(a/b-a'/b') > v(a/b - \alpha_j)$ that $v(a'/b' - \alpha_j)=v(a/b - \alpha_j)$, and our claim is proved.
This claim contradicts the fact that
 $v(F(a,b)) =v( F(a',b')) =v(h)$, and the lemma follows.

\if false
The reader will check  that this inequality implies that 
$v(a/b - \alpha_j ) \geq  v(a'/b' - \alpha_j)$ for all $j$.
This contradicts the fact that
 $v(F(a,b)) =v( F(a',b')) =v(h)$, and our lemma is proved.
\fi

\begin{example}
The above lemma is not correct if the hypothesis that both
$(a,b) $ and $(a',b')$ be primitive solutions is dropped. Indeed,
consider $F(x,y) = (x-y)(x-p^2y)(x-(p^2-p+1)y)^d \in {\mathbb Z}[x,y]$ 
with $(a,b) = (p^2+1,1)$ and $(a',b') = (p,0)$.
\end{example}

Slightly more can be said about the closure of ${\mathcal P}$ in 
${\mathcal Y}/{\mathcal O}_L$. Consider the following two schemes,
${\mathcal U}(\alpha_i)$ attached to a primitive point $(a,b)$ with associated
valuation $t$, and ${\mathcal U}(\alpha_j)$ attached to a primitive point $(a',b')$ with associated
valuation $t'$. We claim that if $v(\alpha_i - \alpha_j) \geq {\rm min}(t,t')$,
then the images of ${\mathcal U}(\alpha_i)$ and ${\mathcal U}(\alpha_j)$
in ${\mathcal Y}$ are equal.
Assume $t' \leq t$.  It follows that $v(a'-\alpha_ib') \geq t'$.
Thus, $v(a'-\alpha_ib') = t'$, and Lemma \ref{decomp} shows that $t=t'$.
We may then define an isomorphism from ${\mathcal U}(\alpha_i)$ to ${\mathcal U}(\alpha_j)$ on the level of rings
$$ {\mathcal O}_L[u',y]/(F'_{m'}(u',y)-\mu) \longrightarrow {\mathcal O}_L[u,y]/(F_m(u,y)-\mu) $$
by setting $u' \mapsto u+ \pi^{-t}(\alpha_i-\alpha_j)y$ and
$y \mapsto y $.

We have thus shown that there exist at most $n$ disjoint disks in
 ${\mathcal O}_L$,
each centered at a root of $F(x,1)$, such that if $\alpha_i$ and $\alpha_j$
belong to the same disk (and have primitive solutions attached to them), then the images of  ${\mathcal U}(\alpha_i)$ and ${\mathcal U}(\alpha_j)$
in ${\mathcal Y}$ are equal. Note now that if a disk contains
a single root of $F(x,1)$, say $\alpha_i$, then by construction
the special fiber of ${\mathcal U}(\alpha_i)$ is a  rational curve,
given by an equation $uy^{n-1}= \overline{\mu}z^n$. In general,
if a disk contains $r$ distinct roots and $\alpha_i$ is one of them,
then the special fiber of ${\mathcal U}(\alpha_i)$ is given 
 by an equation of the form $f_r(u,y)y^{n-r}= \overline{\mu}$,
where $f_r$ is a homogeneous polynomial of degree $r$, coprime to $y$.
%We shall  use later the fact that if the number of disks is $s-i$
%for some $i \leq s/2$, then at least $s-2i$ disks contain only one root
%of $F(x,1)$. 

\begin{emp} {\bf The quotient construction.}
\end{emp}

Let $X/K$ be a smooth proper geometrically connected curve of genus $g$. 
Let $L/K$ be a cyclic Galois extension with Galois group ${\rm  Gal}(L/K) 
= <\sigma>$.
Let 
${\mathcal Y}/{\mathcal O}_L$ be a normal model of $X_L/L$  
such that ${\rm
Gal}(L/K)$ acts on ${\mathcal Y}$, lifting its natural action on ${\rm Spec}({\mathcal O}_L)$.
An example of such a model ${\mathcal Y}$ is the normalization in $L(X)$
of a normal model over ${\mathcal O}_K$ of $X/K$.
Another example is the minimal regular model ${\mathcal Y}/{\mathcal O}_L$
of $X_L/L$. Indeed, the following is well-known.

\begin{emp} Let ${\mathcal Y}/{\mathcal O}_L$ be the minimal regular model
of $X_L/L$.
The map $\sigma$ induces a canonical morphism $X_L \to X_L$ over the map
$\sigma: \Spec(L) \to \Spec(L)$. Since $X_L$ is the  generic fiber of $\yyy$,
the map $\sigma$ induces a birational proper map $\yyy \longrightarrow \yyy \times _{\text{\rm Spec}(\ooo_L)} \Spec(\ooo_L)$ over $\Spec(\ooo_L)$.
By the universal property of a minimal model (\cite[page 310]{C-S}),
this map extends to a 
morphism from $\yyy$ to  $\yyy \times _{\text{\rm Spec}(\ooo_L)} \Spec(\ooo_L)$ over $\Spec(\ooo_L)$. Since $\yyy$ is reduced and separated, this extension is
 unique.
Hence, there exists then  a unique automorphism
$\tau$ of $\yyy$ making the following diagram commutative :
\[
\ba{ccc}
\yyy &\stackrel{\tau}{\longrightarrow}  & \yyy \\
\downarrow &                 & \downarrow \\
{\rm Spec}(\ooo_L) &\stackrel{\sigma}{\longrightarrow}  & {\rm Spec}(\ooo_L) 
\ea
\]
\end{emp}

\begin{emp}
Let $G=< \! \tau \! >$.
The following fact is standard:
Since $\yyy/\ooo_L$ is projective, the quotient $\zzz=\yyy/ G$
can be constructed in the usual way by gluing together the rings of invariants
of $G$-invariant affine open sets of $\yyy$. The scheme $\zzz/\ooo_K$ is  normal   and,
hence, its singular points are closed points of its special fiber.
We let $f : \yyy \longrightarrow \zzz$ denote the quotient map.
\end{emp}

The normal scheme 
$\zzz$ has quotient singularities.
A desingularization $\nu: {\mathcal X} \rightarrow \zzz$   leads to a regular model ${\mathcal X}/{\mathcal O}_K$ of $X/K$.
Let $K^{nr}$ denote the maximal unramified extension of $K$,
and assume now that $K= K^{nr}$.
When $L/K$ is a tamely ramified field extension,
the quotient singularities of $\zzz$
 are well-understood. We recall their properties below,
closely following Viehweg's article \cite{Vie}. We refer the reader to 
his work for more details.
Though he states at the beginning of his paper that he 
considers  only the equicharacteristic case, his proofs of the facts 
listed below are also correct in the mixed characteristic case. 
%We assume that ${\rm ord}( \overline{\tau}) ={\rm ord}(\sigma) = [L:K]$.
If ${\mathcal X}/{\mathcal O}_K$ is any scheme, let us denote by $\overline{\mathcal X}$
the special fiber 
${\mathcal X}\times_{\text{\rm Spec}
({\mathcal O}_K)}\Spec({\mathcal O}_K/(\pi)) $.

\begin{emp}\label{D-def}(\cite[page 303]{Vie})
Let $\overline{\tau}: \overline{\mathcal Y} \rightarrow \overline{\mathcal Y} $
and $\overline{\tau}^{red}: \overline{\mathcal Y}^{red} \rightarrow \overline{\mathcal Y}^{red} $
be the natural morphisms induced by $\tau$. Then the natural map 
\[ \overline{\mathcal Y}^{red}/\!<\overline{\tau}^{red} \! > \longrightarrow \overline{\mathcal Z}^{red}=\overline{\yyy/\!<\tau\!>}^{red} \]
is an isomorphism of schemes over the residue field.
\end{emp}

%Let $\overline{\mathcal Y}^{red} = \bigcup Y_i$ and 
%$\overline{\mathcal Z}^{red}= \bigcup Z_j$. 

For any irreducible component $Y_i \subset \overline{\mathcal Y}$,
let $$D(Y_i):= \{ \mu \in G \;  |  \; \mu(Y_i)=Y_i \}   \mbox{ \rm and } 
I(Y_i):= \{ \mu \in G \; | \; \mu_{|Y_i}={\rm id} \}.$$

\begin{emp} (\cite[page 303]{Vie})
Let $m_i$ be the multiplicity if $Y_i$ in $\overline{\mathcal Y}$ and let
$Z_j:= f(Y_i)$. The multiplicity of $Z_j$ in $\overline{\mathcal Z}$
is equal to $ m_i \cdot [L:K]/|I_i|$.
\end{emp}

Recall the following terminology.
Let $(C \cdot D)$ denote the intersection number on a regular model $\xxx$
of  two divisors $C$ and $D$. Let us call  {\it chain of rational curves on}  $\xxx$
a divisor $D$ such that
\begin{enumerate}
\item $D = \bigcup_{i=1}^q E_i$, $E_i$ smooth and rational
curve for $i=1,\dots,q$.
\item 
$(E_i \cdot E_{i+1})=1$  for all $i=1,\dots,q-1$ and $(E_i \cdot E_j)=0$
for all $j \neq i+1$. Moreover, $(E_i \cdot E_i) \leq -2$ for all $i$.
Let us call  $E_1$ and $E_q$ the end-components of the chain.
\end{enumerate}

Consider again a normal model $\yyy/{\mathcal O}_L$
with an action of Gal$(L/K)$ lifting the action on Spec$({\mathcal O}_L)$.
Assume  that ${\mathcal U}/\ooo_L$ is a smooth open subset 
of  $\yyy/\ooo_L$ such that ${\mathcal U}
$ is invariant under the action of $G$.
%In particular, $ {\mathcal U}_s = \bigcup U_i$ is a union of smooth
%components. 
Let $\zzz:= {\mathcal U}/G$.

\begin{emp} \label{FactV}
(\cite[section 6]{Vie}
There exists a regular scheme $\xxx/\ooo_K$ and a proper 
birational morphism
$\nu : \xxx \rightarrow \zzz$ such that
$\nu$ induces an isomorphism between $\xxx- \{\nu^{-1}(\zzz_{sing}) \}$
and $\zzz - \{ \zzz_{sing} \}$
 and such that, for any $z \in \zzz_{sing}$, 
$\nu^{-1}(z) $ is a connected chain of rational curves. The point $z$
belongs to an end-component of the chain.
Since
${\mathcal U}$ is smooth, we find that
if $z$ is a singular point of $\zzz$, then $\nu^{-1}(z)$ intersects the rest
of the special fiber $\overline{\mathcal X}$ with normal crossings in exactly one point, say on $E_1 $. (Viehweg states in 8.1.d) on page 306 of \cite{Vie} that the model $\xxx$,
obtained by taking the quotient of ${\mathcal U}$ and then resolving the
singularities, has normal crossings.) Let us call the component $E_q$
the {\it terminal component}
of the chain $\nu^{-1}(z) $. The other end-component of  
the chain $\nu^{-1}(z) $ is attached to an irreducible component of 
$\overline{\mathcal X} \setminus \nu^{-1}(z) $.
\end{emp}

 \begin{emp}
  (\cite[section 6]{Vie}) \label{emp.quotient}
Let $f:{\mathcal U} \to \zzz$ denote the quotient map.
Let $z_1,\dots,z_d$ be the closed points of   $\overline{\mathcal Z}$
that are ramification points of the morphism $\overline{f} : \overline{\mathcal Y} \rightarrow\overline{\mathcal Z}^{red}$.
Then $\{z_1,\dots,z_d\}$ is the set of singular points of $\zzz$.
Moreover, if $\nu: \xxx \rightarrow \zzz$ is the desingularization of $\zzz$
described in \ref{FactV}, then the multiplicity  of 
the terminal component on the  chain
$\nu^{-1}(z_i)$ is equal to the number of closed points in the
fiber $\overline{f}^{-1}(z_i)$.
\end{emp}

We now apply the quotient construction to the case where the model
 $\yyy/\ooo_L$ is smooth.
The scheme $\zzz=\yyy/ \! <\!\tau \!>$ has an irreducible special fiber.
The reduced special fiber $\overline{\mathcal Z}^{red}$ is
obtained as the quotient of $\overline{\mathcal Y}$ by the action of $<\! \overline{\tau} \!>$
and is then a smooth and proper curve. 
The multiplicity of $\overline{\mathcal Z}$ in $\zzz$ equals $[L:K]/I(\overline{\mathcal Y})$. The regular  model $\xxx/\ooo_K$ obtained as the minimal 
desingularization of $\zzz$ is thus very simple.

%\vfill \eject
%\input{Chabauty2}

\begin{emp} {\bf Applications of the method of Chabauty-Coleman} 
\end{emp}

We may now apply the method of Chabauty-Coleman 
 to the case of Thue equations. Let $g:=g(X_{F,h})$.

\begin{prop} \label{pro.1}
 Let $X_{F,h}/{\mathbb Q}$ be such that for some prime $p > n$,
$p\nmid h$ and $p\nmid
d^*(F)$.  
Let $K$ be any number field having an unramified prime $\mathfrak P$ 
of norm $p$. 
Assume that 
%the Chabauty rank with respect to $\mathfrak P$ of $X_{F,h}$ over $K$
%is less than
%$g:=g(X_{F,h})$. 
$\Chab(X_{F,h}, K, {\mathfrak P}) < g$.
Then
\begin{eqnarray*}
|X_{F,h}(K)| &\le& (2g-2) \frac{p-1}{p-2} + |\overline{X}_{F,h}({\mathbb F}_p)|.
%&\le& n(n - 2) + p + 1 + 2g\sqrt{p}.
\end{eqnarray*}
\end{prop}

\noindent{\bf Proof: }  As noted in \ref{lem.tame}, $X_{F,h}/{\mathbb Q}_p$ has good
reduction.  Thus, we can apply \ref{to-use}.

\begin{prop} \label{pro.2} Let $X_{F,h}/{\mathbb Q}$ be such that for some prime $p > n$, $p \mid h$
and $p\nmid d^*(F)$.
Let $s$ denote the number of distinct roots of $F(x,1)$ in
$\overline{\mathbb Q}$.
Let $K$ be any number field having an unramified prime $\mathfrak P$ 
of norm $p$. Assume that 
%the Chabauty rank with respect to $\mathfrak P$ of $X_{F,h}$ over $K$
%is less than
%$g:=g(X_{F,h})$. 
$\Chab(X_{F,h}, K, {\mathfrak P}) < g$.
  Then
$|X_{F,h}(K)| \le (2g-2) \frac{p-1}{p-2} + sp$.
\end{prop}

\noindent{\bf Proof: }  Let $X:= X_{F,h}$. As noted in \ref{lem.tame}, 
$X/{\mathbb Q}_p$ has good reduction
after a
tame extension of ${\mathbb Q}_p$.  Thus, we may apply the quotient 
construction to
describe a regular model of $X/{\mathbb Q}^{nr}_p$ over 
${\mathbb
Z}^{nr}_p$.
Let $L = {\mathbb Q}^{nr}_p(\sqrt[n]{h})$.  The extension $L/{\mathbb
Q}^{nr}_p$ is Galois
with cyclic Galois group. Let 
$\xi_n$ be a primitive 
$n$-th root of unity, and denote   by $\sigma: L \rightarrow L$, with
$\sigma(\sqrt[n]{h}) =
\xi_n\sqrt[n]{h}$, a generator of Gal$(L/K)$.
   The morphism $\sigma$ lifts to a morphism $\sigma:
X_{L}
\rightarrow X_{L}$ by setting
$$
\begin{CD}
L[u,v,w]/(F(u,v) - h w^n) @>\sigma>>L[u,v,w](F(u,v)-hw^n)\\
@AAA @AAA\\
L @>\sigma>> L
\end{CD}
$$
with $\sigma(u) = u$, $\sigma(v) = v$ and $\sigma(w) = w$. Let
$
{\mathcal Y}$ denote the normalization of  ${\rm Proj}({\mathcal O}_L[x,y,z]/(F(x,y)-z^n)).
$
 Then ${\mathcal Y}/{\mathcal O}_L$ is the smooth minimal model 
of $X_L/L$ (see \ref{lem.tame}).
The morphism $\sigma: X_{L} \rightarrow X_{L}$ extends to a morphism
$\sigma: {\mathcal Y}
\rightarrow {\mathcal Y}$ by setting
$$
\begin{CD}
{\mathcal O}_L[x,y,z]/(F(x,y)-z^n) @>\sigma>> {\mathcal O}_L[x,y,z]/(F(x,y) - z^n)\\
@AAA @AAA\\
{\mathcal O}_L @>\sigma>> {\mathcal O}_L
\end{CD}
$$
with $\sigma(x) = x$, $\sigma(y) = y$, and $\sigma(z) = \xi_n z$.
When restricted to the special fiber $\overline{\mathcal Y}$ of 
${\mathcal Y}$, the morphism $\sigma$ becomes an
automorphism
$\overline{\sigma}$ over 
$\overline{\mathbb F}_p$ of $\overline{\mathcal Y}$ which lifts
the standard automorphism of ${\rm Proj}(k[x,y,z]/(\overline{F} - z^n))$.  The
quotient map
$\overline{\mathcal Y}
\rightarrow \overline{\mathcal Y}/\langle \overline{\sigma} \rangle$ 
is  ramified over at most
$s$ points.  It
follows from \ref{emp.quotient}
 that the desingularization ${\mathcal X}$ of ${\mathcal Y}/\langle
\sigma \rangle$ has a
special fiber containing at most $s$ (smooth rational) components of
multiplicity one.

Consider now the  minimal regular model ${\mathcal X}_0/{\mathbb Z}_p$
of $X/{\mathbb Q}_p$. A point in $X({\mathbb Q}_p)$ specializes 
in the special fiber $\overline{\mathcal X}_0/{\mathbb F}_p$
to a smooth point, belonging to a geometrically integral  irreducible component
$C/{\mathbb F}_p$ of multiplicity one. Let
$\tilde{{\mathcal X}_0}:= {\mathcal X}_0 \times_{\text{\rm Spec}({\mathbb Z}_p)}\Spec({\mathbb Z}^{nr}_p)$. Since the self-intersection of $C$ in ${\mathcal X}_0$
equals the self-intersection of $C$ in $\tilde{{\mathcal X}_0}$
(see, e.g., \cite{B-L}, 1.4), 
we find that $C$ cannot be contracted in $\tilde{{\mathcal X}_0}$
and, thus, corresponds to a component in the minimal regular model
$\tilde{{\mathcal X}_{00}}$ of $X/{\mathbb Q}_p^{nr}$. Since there is a natural
morphism
  ${\mathcal X} \to \tilde{{\mathcal X}_{00}}$, our description above of the 
special fiber of ${\mathcal X}$ implies 
 that there are at most $s$ components of $\overline{\mathcal X}_0$
that can contain the reduction of a ${\mathbb Q}_p$-point, and that each 
such component is a smooth rational curve. Moreover, each such component $C$
meets the divisor $\overline{\mathcal X}_0 -C$ in exactly 
one ${\mathbb F}_p$-point.
  Hence, the
number of points in $\overline{\mathcal X}_0$ that can be reductions of
${\mathbb Q}_p$-rational points is at
most
$s p$.

\begin{remark} \label{rem.lateproof} Let $K$ be any field with 
a discrete valuation. If $\pi_K \nmid d^*(F)$ and $\pi_K \mid h$, then 
$X/K$ achieves good
reduction over $L:= K(\sqrt[n]{h})$. 
The quotient construction used in the above proof
can be used to  show that  if $M$ 
does not contain $L$, then  
$X_M/M$ does not have good reduction at any maximal ideal of ${\mathcal O}_M$. 
\end{remark}

Let $K$ be any number field, and let ${\mathfrak P}$ be a maximal ideal
of ${\mathcal O}_K$. Let $N(F,h, K, {\mathfrak P})$ denote the number 
of solutions
$(x,y) \in ({\mathcal O}_K)_{\mathfrak P}^2$ of $F(x,y)=h$ with $\gcd(x,y)=1$.

 \begin{prop}  \label{pro.3}
Let $X_{F,h}/{\mathbb Q}$ be such that for some prime $p > n$, $p
\nmid h$ but $p \mid d^*(F)$.
Let $K$ be any number field having an unramified prime $\mathfrak P$ 
of norm $p$. 
Assume that 
%the Chabauty rank with respect to $\mathfrak P$ of $X_{F,h}$ over $K$
%is less than
%$g:=g(X_{F,h})$. 
$\Chab(X_{F,h}, K, {\mathfrak P}) < g$.
 Let $a(p)$
denote the number of
${\mathbb F}_p$-rational points of the affine curve $F(x, y) - h = 0\ {\rm
mod}\ p$.  Then
$$N(F,h,K,{\mathfrak P})
\le (2g-2) \frac{p-1}{p-2} + a(p).
$$
\end{prop}

\noindent{\bf Proof: }  Consider the model ${\mathcal  C}/{\mathbb Z}_p$ given by
$$
{\mathcal  C} = {\rm Proj}({\mathbb Z}_p[x,y,z]/(F - hz^n)).
$$
The special fiber $\overline{\mathcal C}/{\mathbb F}_p$ is a plane projective curve with possible
singularities only at
points $(x: y: z)$ with $z = 0$.  None of the singular points of $\overline{\mathcal
C}$ can be the
reduction of a primitive integral solution of $X_{F,h}({\mathbb Q})$.  Resolve the
singularities of
${\mathcal  C}$ to obtain a regular model ${\mathcal X}/{\mathbb Z}_p$ of $X_{F,h}/{\mathbb
Q}$.  The only points in
$\overline{\mathcal X}/{\mathbb F}_p$ that can be reduction of primitive
points in $X_{F,h}({\mathbb Q})$ are the points in $\overline{\mathcal X}({\mathbb F}_p)$ that correspond to
${\mathbb F}_p$-rational points of $\overline{\mathcal  C}$ with $z \neq 0$.
Apply \ref{to-use}. 

\begin{remark}  The integer $a(p)$ can be bounded using the Weil bound for
singular curves
(where the genus is replaced by the arithmetic genus) and applying it to
each irreducible
components of $\overline{\mathcal  C}$.  When $p$ is not too large compared to $n$
a better bound for
$a(p)$ can be obtained as follows.  Project the curve $\overline{\mathcal  C}/{\mathbb
F}_p$ onto an
${\mathbb F}_p$-rational projective line, using when possible
one of the points at $\infty$
of $\overline{\mathcal  C}({\mathbb F}_p)$.  Then
the projection map has degree at most $n$ (or $n  - 1$
when $\overline{\mathcal  C}({\mathbb F}_p) \neq \emptyset$).
It follows that $a(p) \le np$.
\end{remark}

\begin{prop} \label{pro.4}  Let $X_{F,h}/{\mathbb Q}$ be such that for some prime $p > n$,
$p\mid d^*(F)$ and $p \mid h$.
%Assume that $F(x, 1)$ has distinct roots in ${\mathbb Q}[x]$.  
Let $K$ be any number field having an unramified prime $\mathfrak P$ 
of norm $p$. Assume that 
%the Chabauty rank with respect to $\mathfrak P$ of $X_{F,h}$ over $K$
%is less than
%$g:=g(X_{F,h})$. 
$\Chab(X_{F,h}, K, {\mathfrak P}) < g$.
 Then
$$
N(F,h,K, {\mathfrak P}) \le (2g-2) \frac{p-1}{p-2} + snp.
$$
\end{prop}

\noindent{\bf Proof: }  Let $X:= X_{F,h}$.
Let $L/{\mathbb Q}^{nr}_p$ denote the splitting field
of $F(x, 1)$ over ${\mathbb Q}^{nr}_p$.  Since
$p > n$, the extension $L/{\mathbb Q}^{nr}_p$ is tame and, thus, cyclic.  
Let 
${\mathcal Y}/{\mathcal O}_L$ be the normalization of
$$
{\mathcal  C}/{\mathcal O}_L := {\rm Proj}({\mathcal O}_L[x,y,z]/(F(x,y,z) - hz^n)).
$$
Let $\langle \sigma \rangle = {\rm Gal}(L/K)$.  The morphism
$\sigma$ induces obvious automorphisms
$$
\begin{CD}
{\mathcal Y} @>{\sigma}>> {\mathcal Y}\\
@VVV @VVV\\
{\mathcal C} @>{\sigma}>> {\mathcal C}\\
@VVV @VVV\\
{\rm Spec}({\mathcal O}_L) @>{\sigma}>> {\rm Spec}({\mathcal O}_L).
\end{CD}
$$
We shall denote by $G:= \langle \sigma \rangle$ the group of automorphisms of ${\mathcal Y}$,
resp.\ ${\mathcal C}$, generated by $\sigma$.
Fix a root $\alpha_i$ of $F(x, 1)$ such that there exists a primitive
solution $P=(a, b)$
with $t:= v_L(a - b\alpha_i) = \max_j(v_L(a - b\alpha_j))$.  Let ${\mathcal
U}(\alpha_i)$ denote the
regular open subset of ${\mathcal Y}/{\mathcal O}_L$ described in \ref{emp.ualpha}. Recall that
${\mathcal  U}(\alpha_i) = {\rm
Spec}({\mathcal O}_L[u, y]/(F_m(u,y) - \mu))$. Thus its special fiber
may not be irreducible. The following lemma, whose proof is omitted, describes
the possible components of  $\overline{{\mathcal  U}(\alpha_i)} $.

\begin{lem}  \label{lem.poly}
Let $k$ be any algebraically closed field.  Let
$n \in {\mathbb N}$ with ${\rm char}(k) \nmid n$. 
Let $\xi_n$ denote a primitive $n$-th root of unity in $k$.
 Let $f(x,
y)$ be homogeneous of degree $n$ in $k[x, y]$, and let $\mu
\in k^*$.  Then $f(x, y) - \mu z^n$ factors in $k[x, y, z]$
if and only if there exist $d \mid n$ and $g \in k[x,y]$ with $f
= g^{n/d}$.  Then $f - \mu z^n = \prod^{n/d}_{i=1}(g -
\xi^{id}_{n} \sqrt[n/d]{\mu} z^d)$.
\end{lem}

We will also need the following lemma describing the action of $G$ on
the components of $\overline{{\mathcal  U}(\alpha_i)} $.  Recall the
definitions of $D(Y_i)$ and $I(Y_i)$ in \ref{D-def}.

\begin{lemma} \label{lem.tech}
Let $Y_1, \dots, Y_{n/d}$ denote the irreducible components
of $\overline{\mathcal Y}$ whose generic points belong to ${\mathcal
U}(\alpha_i)$.
Then
$D(Y_{\ell}) = D(Y_j) = G$ and $I(Y_{\ell}) = I(Y_j)$ for all
$\ell, j \in \{1, \dots, n/d\}$.   
\end{lemma}

\noindent {\rm Proof}: Since $p \nmid n$,
the group  of $n$-th roots of unity is contained in ${\mathbb Q}_p^{nr}$,
and acts
on ${\mathcal C}/{\mathcal O}_L$ as follows. A generator $\xi_n$  induces an automorphism $\varphi: {\mathcal
C} \rightarrow {\mathcal C}$ given by:
$$
{\mathcal O}_L[x, y, z]/(F(x, y)- hz^n)
\stackrel{\varphi^*}{\longrightarrow} {\mathcal O}_L[x, y,
z]/(F(x, y) - hz^n),
$$
where $x \mapsto x$, $y\mapsto y$, and $z \mapsto \xi_n
z$.  The automorphism $\varphi$ induces an automorphism
$\varphi: {\mathcal Y} \rightarrow {\mathcal Y}$.  The generator
$\xi_n$ also induces an automorphism $\varphi: {\mathcal
U}(\alpha_i)
\rightarrow {\mathcal U}(\alpha_i)$ given by
$$
{\mathcal O}_L[u, y]/(F_m(u, y) - \mu)
\stackrel{\varphi^*}{\longrightarrow} {\mathcal O}_L[u,
y]/(F_m(u, y) - \mu),
$$
where $u \mapsto \xi^{-1}_n u$ and $y \mapsto \xi^{-1}_n
y$.    
Let
$$
\psi^*: {\mathcal O}_L[x, y]/(F(x,y) - h) \longrightarrow {\mathcal
O}_L[u, y]/(F_m(u,y) - \mu).
$$
be given by $x \mapsto \pi^t_L u + \alpha_i y$, and $y
\mapsto y$.  The induced morphism $\psi: {\mathcal U}(\alpha_i)
\rightarrow {\mathcal C}$ was shown to induce an open immersion
$\psi: {\mathcal U}(\alpha_i) \rightarrow {\mathcal Y}$ in \ref{lem.open}.  The
reader will easily verify that the diagram
$$
\begin{CD}
{\mathcal Y} @>{\varphi}>> {\mathcal Y}\\
@AA{\psi}A @AA{\psi}A\\
{\mathcal U}(\alpha_i) @>{\varphi}>> {\mathcal U}(\alpha_i)
\end{CD}
$$
is commutative.  As above,
let $\sigma $ also denote the automorphisms induced on $\yyy$ and ${\mathcal C}$
by a generator $\sigma$ of ${\rm Gal}(L/K)$.
Since $\sigma$ and $\varphi$ commute, and since $\varphi$
acts transitively on $\{Y_1, \dots, Y_{n/d}\}$, we find that
$D(Y_{\ell}) = D(Y_j)$ and $I(Y_{\ell}) = I(Y_j)$ for all
$\ell, j \in \{1, \dots, n/d\}$.  We let $D := D(Y_j)$ and
$I:= I(Y_j)$.
Note now that $D = G$.  Indeed, if $P= (a,b)$
reduces to $Y_j$ for some $j$, then $\sigma(P)$ reduces to
$\sigma(Y_j)$.  Since $(a, b) \in ({\mathbb Z}_p^{nr})^2$, we find that $P$
reduces to a point in $Y_j \cap \sigma(Y_j)$. Since $P$
reduces to a non-singular point of ${\mathcal U}(\alpha_i)$, we
find
that $Y_j = \sigma(Y_j)$.
This concludes the proof of Lemma \ref{lem.tech}.

\medskip
Consider the following $G$-invariant subset of ${\mathcal Y}$:
$$
{\mathcal V}(\alpha_i) := \bigcap_{\tau \in G} \tau({\mathcal
U}(\alpha_i)).
$$
Let $\tilde{P}$ denote the closure of $P \in X(L)$ in ${\mathcal Y}$.
Then $\tilde{P} \in {\mathcal
U}(\alpha_i)$. Since $P$ is fixed by $\tau$,  $\tau({\mathcal
U}(\alpha_i))$ contains $\tilde{P}$ and, thus, $\tilde{P} \in {\mathcal V}(\alpha_i)$.
 
We need to understand the desingularization of ${\mathcal
V}(\alpha_i)/G$.  Consider first the case where $I = G$. (This case happens for instance if $\alpha_i \in {\mathbb Q}_p^{nr}$.) Then
$$
\overline{{\mathcal V}(\alpha_i)} \longrightarrow \overline{{\mathcal V}(\alpha_i)/G}
$$
is an isomorphism.  Let
\begin{eqnarray}
{\mathcal C} &:=& {\rm Proj} ({\mathcal O}_L[x,y,z]/(F - hz^n))\nonumber\\
{\mathcal D} &:=& {\rm Proj}( {\mathbb Z}^{nr}_p[x,y,z]/(F -
hz^n))\nonumber\\
 {\mathcal D}' &:=& {\rm Proj}({\mathbb Z}_p[x,y,z]/(F - hz^n)).\nonumber
\end{eqnarray}
Let ${\mathcal Z}'/{\mathbb Z}_p$ and  ${\mathcal Z}/{\mathbb Z}^{nr}_p$
denote  the normalization of ${\mathcal
D}'$ and ${\mathcal D}$, respectively.  Clearly ${\mathcal Z} = {\mathcal Y}/G$.  We have the
following commutative diagram:
$$
\begin{CD}
Y_j @.@. \ \ \ \subset @.  \ \ \ \ \ {\mathcal Y} @>{\rho}>>
{\mathcal C} @>>> {\rm Spec}({\mathcal O}_L) \\
@VVV   @.  @. \ \ \ \ \  @VVV  @VVV @VVV\\
Z_{\ell} @. @. \ \ \ \ \subset @.\ \ \ \ \ {\mathcal Z}
@>{\varepsilon}>> {\mathcal D} @>>> 
{\rm Spec}({{\mathbb Z}_p^{nr}})\\
@VVV   @.  @. \ \ \ \ \  @VVV  @VVV @VVV\\
Z_{\ell}' @. @. \ \ \ \ \subset @.\ \ \ \ \ {\mathcal Z}'
@>{\varepsilon'}>> {\mathcal D}' @>>> {\rm Spec}({\mathbb Z}_p).
\end{CD}
$$
The map $\rho$ induces a morphism $\rho_j: Y_j \rightarrow
\rho(Y_j)$, given in coordinates 
by the bottom horizontal map below:
$$
\begin{CD}
{\mathcal O}_L[x, y]/(F(x,y) - h) @>>> 
{\mathcal O}_L[u, y]/(F_m(u,y) -\mu)\\
@VVV @VVV \\
{\mathcal O}_L[x, y]/(\pi_L, x- \alpha_iy) @>>> 
{\mathcal O}_L[u, y]/(\pi_L, F_m(u,y) -\mu).\\
\end{CD}
$$
This morphism is clearly of degree at most $n$.
The  morphism
$$
\varepsilon'_{\ell}: Z'_{\ell} \rightarrow
\varepsilon'_{\ell}(Z'_{\ell})
$$
induced by  $\rho_j$ is
also of degree at most $n$.  The curve % 
$\varepsilon'_{\ell}(Z'_{\ell})/{\mathbb F}_p$ is a smooth projective
line. The primitive integral point $(a,b)$ cannot reduce to the intersection point $Q$ of all the components of the special fiber of  ${\mathcal D}'$.
The morphism $\varepsilon'_{\ell}$ 
is defined over ${\mathbb F}_p$, and there are at most $np$ ${\mathbb F}_p$-rational 
points in the preimage of $\varepsilon'_{\ell}(Z'_{\ell}) \setminus \{ Q \}$.
We conclude that at most $np$
points in the image of ${\mathcal V}(\alpha_i)$ in ${\mathcal Z}'$ can be
residue classes of primitive integral points.

Let us now consider the case where $I  \subsetneq
D$.  Then the image $Z_\ell$ of $Y_j$ in ${\mathcal Z} = {\mathcal Y}/G$ has
multiplicity $|D|/|I| >1$, and \ref{emp.quotient} indicates that to count the 
components of multiplicity one (in a desingularization of ${\mathcal Z}' $)
which contain the reduction of primitive integral points, one first needs
 to count the number of
totally ramified points in the branch locus of  $Y_j \rightarrow Z_{\ell}$.
Consider
the diagram
$$
\begin{CD}
{\mathcal Y} @>{\sigma}>> {\mathcal Y}\\
@AA{\psi}A @AA{\psi}A\\
{\mathcal V}(\alpha_i) @>{\sigma}>> {\mathcal V}(\alpha_i),
\end{CD}
$$
where $\sigma: {\mathcal V}(\alpha_i) \to {\mathcal V}(\alpha_i)$ is  
defined so that the diagram commutes.
Consider an open set ${\mathcal U}$ of ${\mathcal V}(\alpha_i)$ that is dense in each fiber
and is a special open set of ${\mathcal U}(\alpha_i)$.
We find that on the level of rings, $\sigma: {\mathcal U} \to {\mathcal U}(\alpha_i)$ induces the top horizontal  map below
$$
\begin{CD}
{\mathcal
O}_L[u, y]/(F_m(u,y) - \mu) @>{\sigma}>> S^{-1}({\mathcal
O}_L[u, y]/(F_m(u,y) - \mu))\\
@AAA @AAA\\
{\mathcal
O}_L[x, y]/(F(x,y) - h) @>{\sigma}>>  {\mathcal
O}_L[x, y]/(F(x,y) - h)\\
@AAA @AAA\\
{\mathcal O}_L @>{\sigma}>> {\mathcal O}_L,
\end{CD}
$$
with $\sigma(x) = x$ and $\sigma(y)=y$
so that  $\sigma(\pi_L^{t}u+\alpha_iy) = \pi_L^{t}u+\alpha_iy$.
Since $\sigma(\pi_L^{t}u+\alpha_iy) = \sigma(\pi_L^{t})\sigma(u)+
\sigma(\alpha_i)y$, we find that 
$$\sigma(u) = \frac{\pi_L^{t}}{\sigma(\pi_L^{t})} u + \frac{\alpha_i -
  \sigma(\alpha_i)}{\sigma(\pi_L^{t})} y.$$
(Note that both
$\pi_L^{t}/\sigma(\pi_L^{t})$ and $(\alpha_i -
\sigma(\alpha_i))/\sigma(\pi_L^{t}) $ belong to ${\mathcal O}_L$.)  By
hypothesis, $\overline{\sigma} := \sigma_{|Y_\ell}$ does not act trivially on $Y_{\ell}$.  The points
where the morphism $Y_{\ell} \to Y_{\ell}/<\overline{\sigma}>$ is totally
ramified is the set of fixed points of the map $\overline{\sigma}$.
On the plane curve $\overline{F_m(u,y) - \mu}=0$, the automorphism
$\overline{\sigma}$ is given by $u\mapsto cu+dy$ and $y\mapsto y$,
for some $c,d \in k$. Thus the fixed points of
$\overline{\sigma}$ lie on the line $(c-1)u+dy=0$, and we
 find that there are at most $n$ such points.
%%, since the intersection
%%of a line with a plane curve of degree $n$ contains at most $n$
%%distinct points. 
 Let now $\nu : {\mathcal X} \to {\mathcal Z}$ denote
the minimal desingularization of ${\mathcal Z}$. As we recalled in
\ref{emp.quotient}, the special fiber of ${\mathcal X}$ contains at
most $n$ components of multiplicity one, each smooth and rational, and
each meeting the rest of the special fiber in a single point.

Consider now the minimal regular model ${\mathcal X}_0/{\mathbb Z}_p$
of $X/{\mathbb Q}_p$. A point in $X({\mathbb Q}_p)$ specializes in the
special fiber $\overline{\mathcal X}_0/{\mathbb F}_p$ to a smooth
point, belonging to a geometrically integral irreducible component
$C/{\mathbb F}_p$ of multiplicity one. Let $\tilde{{\mathcal X}_0}:=
{\mathcal X}_0 \times_{\text{\rm Spec}({\mathbb Z}_p)}\Spec({\mathbb
  Z}^{nr}_p)$.
 Since the self-intersection of such a component $C$ in ${\mathcal X}_0$
equals the self-intersection of $C$ in $\tilde{{\mathcal X}_0}$
(see, e.g., \cite{B-L}, 1.4), we
find that $C$ cannot be contracted in $\tilde{{\mathcal X}_0}$ and,
thus, corresponds to a component in the minimal regular model
$\tilde{{\mathcal X}_{00}}$ of $X/{\mathbb Q}_p^{nr}$. Since there is
a natural morphism ${\mathcal X} \to \tilde{{\mathcal X}_{00}}$, our
description above of the special fiber of ${\mathcal X}$ implies that
there are at most $n$ components of $\overline{\mathcal X}_0$ that can
contain the reduction of a ${\mathbb Q}_p$-point, and that each such
component is a smooth rational curve. Moreover, each such component
$C$ meets the divisor $\overline{\mathcal X}_0 -C$ in exactly one
${\mathbb F}_p$-point.  Hence, the number of points in
$\overline{\mathcal X}_0$ that can be reductions of ${\mathbb
  Q}_p$-rational points is at most $n p$.

Since the contribution of an open set of the form ${\mathcal V}(\alpha_i)$
to the number of reductions of  primitive integral points   
in the special fiber of the model ${\mathcal X}_0$ is bounded by $np$, and since
the primitive integral points are contained in at most $s$ 
such open sets (\ref{lem.open}),
we find that   the reduction of the primitive integral points   
in the special fiber of the model ${\mathcal X}_0$ consists in at
most $snp$ points. 
This concludes the proof of Proposition \ref{pro.4}.

\medskip
Let us now state our main theorem.   Let $N(F,h)$ denote the number of solutions
$(x,y) \in {\mathbb Z}^2$ of $F(x,y)=h$ with $\gcd(x,y)=1$.

\begin{thm}  \label{thm.main}
Let $p$ be a prime  with $n < p < 2n$. 
 Assume that
the Chabauty rank with respect to $(p)$ of 
$  X_{F,h}/{\mathbb Q} $ is less than $ g:= g(X_{F,h})$.  Then
$$
N(F,h) \le 2n^3 - 2n -3.
$$
More precisely, choose a prime $p$ with $n < p < 2n$.
\begin{enumerate}
\item[a)]  If $p \nmid h$ and $p \nmid d^*(F)$, then $|X_{F,h}({\mathbb Q})|
\le 2g + s-5 +2n(n-1)$.
\item[b)]  If $p \mid h$ and $p \nmid d^*(F)$, then $|X_{F,h}({\mathbb Q})|
\le 2g -5 + 2sn$.
\item[c)]  If $p \nmid h$ and $p \mid d^*(F)$, then $N(F,h,{\mathbb Q},p)
\le 2g + s-5 +n(2n-1) $.
\item[d)]  %Assume that $F(x, 1)$ does not have multiple
%roots. 
If $p \mid h$ and $p \mid d^*(F)$, then $N(F,h,{\mathbb Q},p) \le
2g + s-5 +sn(2n-1)$.
\end{enumerate}
In particular, if the Mordell-Weil rank of 
$  X_{F,h}/{\mathbb Q} $ is less than $ g$, then 
$N(F,h) \le 2n^3 - 2n -3$.
\end{thm}

\noindent{\bf Proof:} We apply our previous results using the estimate
$p \le 2n - 1$ and $s \leq n$. 
The term $(2g-2)(p-1)/(p-2)$ is bounded by $2g+s-5$.
To prove a), apply \ref{pro.1},
and bound $|\overline{X}_{F,h}({\mathbb F}_p)|$ using a projection from 
a smooth point of  $\overline{X}_{F,h}({\mathbb F}_p)$ to a ${\mathbb F}_p$-line; we find in this case that $|\overline{X}_{F,h}({\mathbb F}_p)| \leq (n-1)(p+1)$. (Note that we may assume that 
$\overline{X}_{F,h}({\mathbb F}_p)$ is not empty because otherwise $X_{F,h}({\mathbb Q}) $ is empty.)
To prove c), apply \ref{pro.3},
and bound $a(p)$ using a projection from a point of  $\overline{X}_{F,h}({\mathbb F}_p)$ 
to a line; we find again that $a(p) \leq (n-1)(p+1)$.
To prove b) and d), use \ref{pro.2}  and \ref{pro.4}.

%% Tom
Note now that by Bertrand's postulate, there exists a prime $p$ with 
$n < p < 2n$. If the Mordell-Weil rank of 
$  X_{F,h}/{\mathbb Q} $ is less than $ g$, then 
the Chabauty rank with respect to $(p)$ of 
$  X_{F,h}/{\mathbb Q} $ is also less than $ g$, and we find that $N(F,h) \le 2n^3 - 2n -3$.

\begin{remark} \label{emp.othersets} 
Recall that for any $a \in {\mathbb Z}$,
the curve $X_{F,hp^{an}}$ is isomorphic over ${\mathbb Q}$
to $X_{F,h}$. It follows that the Chabauty rank of $X_{F,hp^{an}}/{\mathbb Q}$
is equal to the Chabauty rank of $X_{F,h}/{\mathbb Q}$. Thus, when
the method of Chabauty-Coleman can be applied to bound $N(F,h,{\mathbb Q},p)$,
it can also be used to bound the size of subsets of $X_{F,h}({\mathbb Q})$
other than the subset of primitive $p$-integral solutions. For instance, consider the case where $F(x,1)$ has distinct roots and let, for $i \geq 0$, 
$$S_i:= \{  (x:y:z) \in X_{F,h}({\mathbb Q}) \mid x,y,z \in {\mathbb Z}_{(p)}, p \nmid z, (x,y) = p^i{\mathbb Z}_{(p)} \}.
$$
The set $S_0$ is the set of primitive $p$-integral points of $X_{F,h}/{\mathbb Q}$.
 The set $S_i$ is in bijection with the set of primitive
$p$-integral points of $X_{F,hp^{-in}}/{\mathbb Q}$.
The set $$T_i:= \{  (x:y:z) \in X_{F,h}({\mathbb Q}) \mid x,y,z \in {\mathbb Z}_{(p)},  (x,y, z) = {\mathbb Z}_{(p)}, v_p(z) =i  \}
$$
 is in bijection with the set of primitive
$p$-integral points on $X_{F,hp^{in}}/{\mathbb Q}$.
\end{remark}

\if false
Recall that the map $\rho: X \rightarrow
C_{F,h}$ is an isomorphism outside the points at infinity of
$C_{F,h}$.  Let us call $P \in X(K)$ a primitive integral point
if $\rho(P) = (x: y :1)$ with $x, y \in {\mathcal O}_K$
and $x {\mathcal O}_K + y{\mathcal O}_K = {\mathcal O}_K$.
When $K$ has a valuation $v$, we define more generally the
subset $X_{F,h}(i)(K) \subseteq X_{F,h}(K)$
for $i \in {\mathbb Z} \setminus\{0\}$ as follows.
A $K$-rational point $P $ belongs to $ X_{F,h}(i)(K)$ if
 $\rho(P) = (x:y:z)$  with
$ x, y, z \in {\mathcal O}_K$, $x{\mathcal O}_K + y{\mathcal O}_K +
z{\mathcal O}_K = {\mathcal O}_K$,  and  $v(z) = -i $
 when $ i < 0$. When $i \geq 0$, $P$ belong to  $X_{F,h}(i)(K)$
if $v(z) = 0$  and $ x{\mathcal O}_K + y{\mathcal
O}_K = \pi^i {\mathcal O}_K $.
We show in \ref{emp.othersets} that each set $X_{F,h}(i)$ has order bounded by a
constant depending on $n$ only, provided that the method of
Chabauty-Coleman can be applied.
\fi

\begin{remark} \label{rem.toomanyP1}
The method of Chabauty-Coleman bounds the number of rational points in
$X_{F,h}({\mathbb Q})$ in terms of
the number $N$ of irreducible components
(of the special
fiber of a regular
model) of multiplicity one 
that contain ${\mathbb F}_p$-points.  Thus, the method provides a bound
for $|X_{F,h}({\mathbb Q})|$
in terms of $n$ if $N$ can be bounded in terms of $n$ only.  As the
following example with $n = 6$
shows, the integer $N$ is not bounded in terms of $n$ in general.

Consider a genus 2 curve $Y/{\mathbb Q}$ with a regular model ${\mathcal Y}/{\mathbb
Z}_p$ whose special
fiber consists of a chain of projective lines over ${\mathbb F}_p$ with one elliptic
curve attached to each
end of the chain.  All intersections in $\overline{\mathcal Y}$ are transverse and all
components have
multiplicity 1.  Thus ${\rm Jac}(Y)/{\mathbb Q}$ has good reduction modulo $p$.
  An example of
such a curve when $p\ge 7$ is the curve $[I_0 - I_0 - m]$ in
\cite{N-U} given by
$$
y^2 = (x^3 + \alpha x + 1)(x^3 + \beta p^{4m} x + p^{6m}), \mbox{ with } \alpha,
\beta \in {\mathbb Z}_p^*.
$$
The number of rational curves in the special fiber of the minimal model
${\mathcal Y}/{\mathbb Z}_p$ is
$m$ (and not $m+1$ as stated in \cite{N-U}).  
This can be shown using Liu's algorithm (\cite{Liu}, Thm.\ 1 and Prop.\ 2).
Consider now the
curve $X/{\mathbb Q}$ given by
$$
y^6 = (x^3 + \alpha x + 1) (x^3 + \beta p^{4m} x + p^{6m}).
$$
We claim that the minimal model ${\mathcal X}/{\mathcal O}_K$ has a special fiber
$\overline{\mathcal X}$ consisting
of two curves of genus 4 linked by three chains of rational curves, as in
the picture below on the left. The special fiber $\overline{\mathcal Y}$
is represented below.

\centerline{\BoxedEPSF{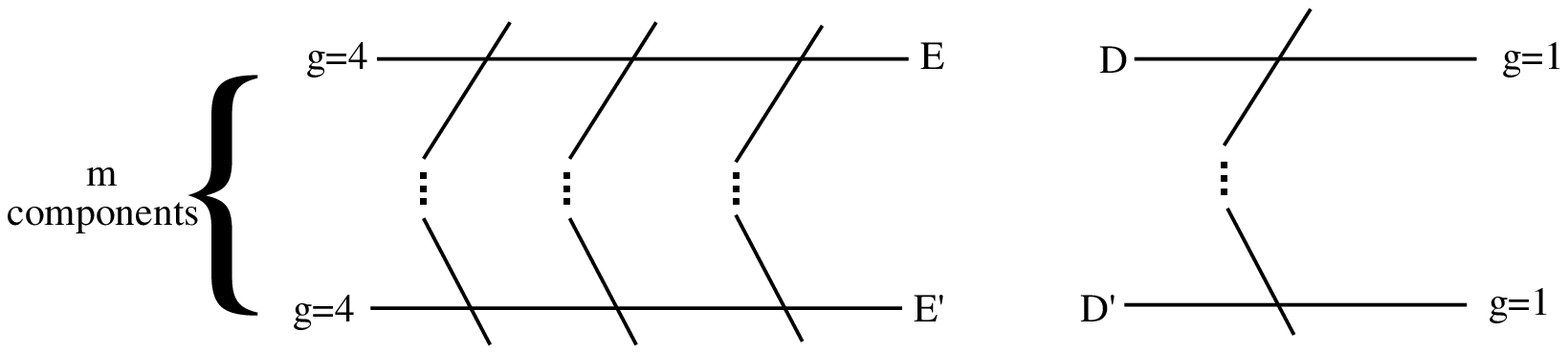  scaled 600}}

Let us only give a brief sketch of the proof of this claim.
First, the map ${\varphi}: E \rightarrow D$ is the natural map between the
curve $E$ given by $y^6
= (x^3 + \overline{\alpha} x + 1)x^3$ and the curve $D$ given by $y^2 =
(x^3 + \overline{\alpha} x +
1)x^3$.  Making the change of variable $u:= x/t^{2m}$, the curve $E'$ is
given by $y^6 = x^3 +
\overline{\beta}x + 1$ and the curve $D'$ by $y^2 = x^3 + \overline{\beta}x
+ 1$.  The map ${\varphi}' : E' \rightarrow D'$ is the natural map.  The maps $\varphi$ and
$\varphi'$ are not ramified
over the points at $\infty$ of $D$ and $D'$.  Thus the graph associated with a
regular model of $X$
having normal crossings must contain at least two independent loops.  Since $g(X) =
10$ and $g(E) = g(E')
= 4$, we conclude that $X$ has semistable reduction over ${\mathcal O}_K$.  The
automorphism $\sigma:
(x, y) \rightarrow (x, \xi_3 y)$ induces an action on the minimal regular
model ${\mathcal X}/{\mathcal
O}_K$ of $X/K$.  It can be shown that ${\mathcal X}/\langle \sigma \rangle$ is
a semi-stable regular model for
$Y/K$.  Since the minimal model ${\mathcal Y}/{\mathcal O}_K$ has a special
fiber with $m$
projective lines, we conclude that $\overline{\mathcal X}$ contains $3m$
projective lines.
\end{remark}

\section{Some refinements} \label{sec.ref}

We use in this section the fact that
any subfield of ${\mathbb Q}(\xi_{p-1})$ has  an unramified
prime of norm $p$ to obtain in some cases bounds for $N(F,h, {\mathbb Q},p)$
of the form $O(n^2)$.

\begin{thm}  \label{thm.ref1}
Let $n  \ge 5$ be a prime so that $p:= a n + 1$ is
also prime for some $a >
1$.  Assume that the Mordell-Weil rank of ${\rm Jac}(X_{F,h})({\mathbb Q}))$ is
less than $(n - 3)/2$. 
Then
\begin{enumerate}
\item[a)]  If $p \nmid h$ and $p \nmid d^*(F)$, then $|X_{F,h}({\mathbb Q})| 
\le (a
+ 2)n - (a + 1)$.
\item[b)] If $p \mid h$ and $p \nmid d^*(F)$, then $|X_{F,h}({\mathbb Q})| 
\le (a +2)n- 2$.
\item[c)] If $p \nmid h$ and $p \mid d^*(F)$, then $N(F,h, {\mathbb Q},p) 
\le (a
+ 1)(n -  1)$.
\item[d)]  
 If $p\mid h$ and $p\mid d^*(F)$, then $N(F,h,{\mathbb Q},p) \le an^2+2n-3 $.
\end{enumerate}
\end{thm}

\noindent{\bf Proof:} Let $X:= X_{F,h}$.
 Recall (\cite{P-S}, 13.4) that when $n$ is prime,
$$
{\rm rank}_{{\mathbb Z}}({\rm Jac}(X/{\mathbb Q}))(n - 1) = {\rm
rank}_{{\mathbb Z}}{\rm
Jac}(X/{\mathbb Q}(\xi_{n})).
$$
Thus, our hypothesis implies that the Mordell-Weil rank of ${\rm
Jac}(X/{\mathbb Q}(\xi_{n}))$
is less than $g(X)$.  We may then apply the results of the previous
section.  We bound $(2g-2)(p-1)/(p-2) $ by $n^2 -2n-3$, and $s$ by $n$.
Let $u$ denote the
number of points in $X({\mathbb Q})$ with $z = 0$. Clearly, $u \leq n$. Let $v:=|X({\mathbb Q})|
- u$.  Then
$$
|X({\mathbb Q}(\xi_{n}))| \ge u +
n v.
$$
For part (a), we use the bound \ref{pro.1}:
\begin{eqnarray*}
|X({\mathbb Q}(\xi_{n}))| \leq  n^2- 2n -3 + (n - 1)(p + 1).
\end{eqnarray*}
(To bound $|\overline{X}({\mathbb F}_p)|$, use a projection from a point in
$\overline{X}({\mathbb F}_p).$).
It follows that
$$
u+ n v \le (a + 1)n^2 - an  - 5.
$$
Hence,
$$
v \le (a + 1) n- a - 5/n - u/n.
$$
Thus, $v \le (a + 1)(n - 1)$. We find that
$$
|X({\mathbb Q})| = u + v \le n+ (a + 1)(n - 1).
$$
To prove part (b), we use \ref{pro.2},
$$
|X({\mathbb Q}(\xi_{n}))| \le n^2- 2n -3 + n p.
$$
Thus,
$$
u + n v \le (a + 1)n^2 - n -3.
$$
Hence, $v \le (a + 1)n - 2$, and $|X({\mathbb Q})| \le (a + 2)n - 2$.
To prove part c), we use \ref{pro.3}:
$$
nN(F,h,{\mathbb Q},p) \leq 
N(F,h,{\mathbb Q}(\xi_{n}),p) \le n^2- 2n -3 + (n - 1)(p + 1).
$$
(To bound $a(p)$, we use the fact that $a(p) \le |\overline{X}({\mathbb
F}_p)|$, and project
$\overline{X}({\mathbb F}_p)$ to a line using one of its points).
Finally, to prove part d), we use \ref{pro.4}:
$$
nN(F,h,{\mathbb Q},p) \leq N(F,h,{\mathbb Q}(\xi_{n}),p) \le n^2- 2n -3 +n^2p.
$$

\begin{thm}  \label{thm.ref2}
 Let $p \ge 5$ be prime and let $n:= p - 1$.  Let $X:= X_{F,h}$.
Assume that $\Chab(X,{\mathbb Q}(\xi_{p-1}), (p)) < g(X)$. This is the case,
for instance, if the Mordell-Weil rank of $X/{\mathbb Q}$ is less than
 $(s-2)/2$.
  Then
\begin{enumerate}
\if false
\item[a)]  If $p \nmid h$ and $p \nmid d^*(F)$, then 
$N(F,h,{\mathbb Q},p) \le 4n-3$ and 
$|X({\mathbb Q})| \le 5n
- 3$.
\item[b)]  If $p \mid h$ and $p \nmid d^*(F)$, then 
$N(F,h,{\mathbb Q},p) \le 4n-3$ and $|X({\mathbb Q})| \le 5n-3$.
\fi
%%Tom  Both lines gave the same bound. (I hope the bounds are correct.)
\item[a)]  If   $p \nmid d^*(F)$, then 
$N(F,h,{\mathbb Q},p) \le 4n-3$ and 
$|X({\mathbb Q})| \le 5n
- 3$.
\item[b)]  If $p \nmid h$ and $p \mid d^*(F)$, then $N(F,h,{\mathbb Q},p) \le 4n-3$.
\item[c)] If $p \mid h$ and $p \mid d^*(F)$, then $N(F,h,{\mathbb Q},p) \le 2n^2 +4n -5$.
\end{enumerate}
\end{thm}

\noindent{\bf Proof:}  We proceed as in the previous theorem,
with $
|X({\mathbb Q}(\xi_{n}))| \ge u +
n v/2$.
For part (a), we use the bound \ref{pro.1}:
$$
|X({\mathbb Q}(\xi_{n}))| \le n^2- 2n -3 + (n - 1)(p + 1).
$$
(To bound $|\overline{X}({\mathbb F}_p)|$, use a projection from a point in
$\overline{X}({\mathbb F}_p).$).
It follows that
$
u+ n v/2 \le 2n^2 -  n  - 5$.
Hence,
$
v \le 4n- 2 - 10/n - 2u/n$.
Thus, $v \le 4n-3$. We find that
$$
|X({\mathbb Q})| = u + v \le n+ 4n - 3.
$$
To prove part (b), we use \ref{pro.2},
$$
|X({\mathbb Q}(\xi_{n}))| \le n^2- 2n -3 + n p.
$$
Thus, $
u + n v/2 \le  2n^2 - n -3$.
Hence, $v \le 4n-3$, and $|X({\mathbb Q})| \le 5n-3$.
To prove part c), we use \ref{pro.3}:
$$
nN(F,h,{\mathbb Q},p)/2 \leq N(F,h,{\mathbb Q}(\xi_{n}),p) \le n^2- 2n -3 + (n - 1)(p + 1).
$$
(To bound $a(p)$, we use the fact that $a(p) \le |\overline{X}({\mathbb
F}_p)|$, and project
$\overline{X}({\mathbb F}_p)$ to a line using one of its points).
Finally, to prove part d), we use \ref{pro.4}
$$
nN(F,h,{\mathbb Q},p)/2 \leq N(F,h,{\mathbb Q}(\xi_{n}),p) \le n^2- 2n -3 +n^2p.
$$

The assertion that, if the Mordell-Weil rank of
$X/{\mathbb Q}$ is less than $(s-2)/2$, then
the Chabauty rank of ${\rm Jac}(X/{\mathbb
  Q}(\xi_{p-1}))$ is less than $g(X)$, is a consequence of the
following general fact proved in \ref{pro.dim} below.  Let $X/{\mathbb
  Q}$ denote the smooth proper model of the affine curve given by an
equation $hy^n = f(x)$, with $n \mid \mbox{\rm deg}(f)$.  Write $f(x)
= \prod^s_{i=1}(x - a_i)^{n_i}$ with $\prod_{i\ne j} (a_i - a_j) \ne
0$, and assume that $\gcd(n,n_i) < n$ for all $i$.  Lemma 13.4 in
\cite{P-S} states that when $n$ is prime, then
$$
{\rm rank}_{\mathbb Z} ({\rm Jac}(X/{\mathbb Q}))(n - 1) = 
{\rm rank}_{{\mathbb Z}} ({\rm
Jac}(X/{\mathbb Q}(\xi_n))).
$$
Fix $\xi_n$, a primitive $n$-th root of unity. Denote by $\sigma$
the automorphism of $X/{\mathbb Q}(\xi_n)$ induced by $(x, y) \mapsto
(x, \xi_n y)$.  The proof of 13.4 when $n$ is prime relies strongly on
the fact that the minimal polynomial of $\sigma$ acting on ${\rm
  Jac}(X)({\mathbb Q}(\xi_n))$ is the minimal polynomial of $\xi_n$.
When $n$ is not prime, the minimal polynomial of $\sigma$ on ${\rm
  Jac}(X)({\mathbb Q}(\xi_n))$ divides $t^{n-1} + t^{n-2} + \cdots + t
+ 1$ and may not be irreducible; the proof of 13.4 does not apply.  We
can nevertheless show that the following weaker statement holds.

\begin{proposition} \label{pro.dim} Assume that the Chabauty rank (with respect to any prime)
of $X/{\mathbb Q}(\xi_n)$
is 
equal to $g(X)$.  Then the Mordell-Weil rank of $X/{\mathbb Q}$ is at least
equal to $(s- 2)/2$.

\end{proposition}

\noindent {\rm Proof}: 
Let $\phi(t) = (t^n - 1)/(t - 1)$.  Let $\varphi_d(t)$ denote the
$d$-th cyclotomic polynomial, so that $\phi(t) =
\prod_{\stackrel{d|n}{d\ne 1}} \varphi_d(t)$.  Let $\phi_d(t) :=
\phi(t)/\varphi_d(t)$.  Let $d \mid n$, and consider the abelian
variety
$$
A_d := {\rm Im}(\phi_d(\sigma)) \subset {\rm Jac}(X)/{\mathbb Q}(\xi_n).
$$
(Note that when $d < \gcd(n,n_i)$ for some $i$, it may happen that
$A_d$ is trivial.)  It is clear that $A_d \subset {\rm
  Ker}(\varphi_d(\sigma))$.  If $d \ne d'$, $\varphi_d(t)$ and
$\varphi_{d'}(t)$ are coprime, and we conclude that $A_d \cap A_{d'} =
\{0\}$.  Since the polynomials $\{\phi_d(t)\}_{d \mid n}$ are coprime,
we can find $\{a_d \in {\mathbb Z}, d \mid n, d \ne 1\}$ such that
$\sum a_d \phi_d(t) = 1$.  Hence, given $P \in {\rm Jac}(X)$,
$$
P = \sum \phi_d(\sigma)(a_d P) \in 
\langle A_d, d \mid n, d \ne 1\rangle \subseteq
{\rm
Jac}(X)
$$
and, thus, ${\rm Jac}(X) = \oplus_{\stackrel{d \mid n}{d\ne 1}} A_d$.

We claim now that $A_d$ is an abelian variety defined over ${\mathbb
  Q}$.  Indeed, let $P = (a, b)$, where $a,b \in \Qb$, be a solution
of $y^n = f(x)$.  Let $\mu \in {\rm Gal}(\overline{\mathbb Q}/
{\mathbb Q})$.  Then $\mu(\xi_n) = \xi^c_n$ for some $c \in {\mathbb
  Z}$ with $ (c, n) = 1$.  It follows that on $X(\overline{\mathbb
  Q})$, we have
$$\sigma^c \circ \mu = \mu \circ \sigma$$
for some $c \in \mathbb Z$.
If an element $Z$ of ${\rm Jac}(X)(\overline{\mathbb Q})$ is of the form
$$
Z = \phi_d(\sigma)(\sum^s_{i=1}a_i P_i), \mbox{ \rm with } \ P_i \in 
X({\overline{\mathbb Q}}),
$$
then $\mu(\phi_d(\sigma)(\sum^s_{i=1}a_iP)) =
\phi_d(\sigma^c)(\sum^s_{i=1} a_i \mu(P_i))$.  Since $(c, n) = 1$ and
$d \mid n$, we find that $\phi_d(t)$ divides $\phi_d(t^c)$.  Hence,
$\mu(Z) \in A_d(\overline{\mathbb Q})$, for all $\mu \in {\rm
  Gal}(\overline{\mathbb Q}/{\mathbb Q})$.  It follows that $A_d$ is
defined over ${\mathbb Q}$.  We determine the dimension of $A_d$ below.

\begin{lemma} \label{lem.dim}
Suppose that $d \mid n$ and that $d > \gcd(n,n_i)$ for all $i=1,\dots,s$. Let   $\varphi(d)$ denote the
Euler
$\varphi$-function.
Then $\dim(A_d) = \varphi(d)(s - 2)/2$.
 \end{lemma}

\noindent{\bf Proof:}  By construction, $\sigma_{|A_d}$ is such that
$\varphi_d(\sigma_{|A_d}) = 0$.  The characteristic polynomial ${\rm
  char}(t)$ of $\sigma$ acting on $H_1(X({\mathbb C}), {\mathbb C})$
is computed in \cite{Lor}, 4.1:
$$
{\rm char}(\sigma)(t) = \phi(t)^{s-2}\prod_{i=1}^s \left(
  \frac{x^{\gcd(n,n_i)}-1}{x-1} \right)^{-1}.
$$
Hence, ${\rm rank}_{\mathbb Z}({\rm Ker}(\varphi_d(\sigma)_{|H_1(X({\mathbb C}),
{\mathbb C})}) =
(s - 2)\varphi(d)$.
Using the duality between $H_1(X({\mathbb C}), {\mathbb C})$ and $H^1(X({\mathbb C}),
{\mathbb
C})$ as well as the fact that
$$
0 \rightarrow H^0(X({\mathbb C}), \Omega_X) \rightarrow H^1(X({\mathbb C}), {\mathbb C})
\rightarrow H^1(X, {\mathcal O}_X) \rightarrow 0
$$
is exact, with $H^0(X({\mathbb C}), \Omega_X)$ and $H^1(X, {\mathcal O}_X)$ related by
Serre duality,
we find that $\dim A_d = \varphi(d)(s - 2)/2$.  This concludes the proof of
Lemma \ref{lem.dim}.

\medskip 
Now, the reader will easily check that the proof of Lemma 13.4 in
\cite{P-S} can be used, {\it mutatis mutandis}, to show that, for any
$d \mid n$,
$$
{\rm rank}_{\mathbb Z}(A_d({\mathbb Q}))\varphi(d) = {\rm rank}_{\mathbb Z}(A_d({\mathbb
Q}(\xi_d))).
$$
In particular, if the Chabauty rank of ${\rm Jac}(X)/{\mathbb Q}$   equals
$g(X)$, then
$$
{\rm rank}_{\mathbb Z}(A_n({\mathbb Q}(\xi_n))) \ge \dim(A_n),
$$
so that ${\rm rank}_{\mathbb Z} (A_n({\mathbb Q})) \ge (s - 2)/2$.
This concludes the proof of \ref{pro.dim}.

\begin{remark}
The bound obtained in Theorem \ref{thm.main} 
is probably  not optimal.  The interest in Theorem \ref{thm.ref1} and 
\ref{thm.ref2}
is that they exhibit cases where the bound for $N(F,h)$ is $O(n)$.
There are no known examples in the literature of a family of Thue
equations $(F_n(x,y) = h_n)_{n=1}^{\infty}$ with the degree of $F_n$
going to infinity for which $N(F_n,h_n) > O(\deg F_n)$.
The following simple examples
of Thue equations
with $|N(F, h)| \geq n+1$ are well-known.  Take
$$
F(x, y) = \prod^n_{i=1}(x - a_i y) + hy^n, \mbox{ with } \prod_{i\ne
j}(a_i - a_j) \ne 0,   a_i \in {\mathbb Z}.
$$
Then $\{(a_i, 1) , i = 1, \dots, n\}$ are primitive
solutions (note that $(\ell, 0)$ is also a solution if
$\ell^n = h$).  If $n$ is even, $\{(-a_i, -1), i = 1, \dots,
n\}$ are also solutions. Now pick $q \in {\mathbb Z}$ and set $a_1 := q^{n-1}$. Take $h :=
\prod_{i=2}^n (1-a_iq)$.
Then $(1,q)$ is an additional primitive solution. The determination of all solutions of such a Thue equation
can sometimes be done (see \cite{Heu} and its bibliography list).
 The current 
record on the number of integral points on plane curves can be found in 
\cite{R-V} and its review in the Math Reviews.

%One certainly expects that 
%the points $(a_i,1) - (a_1,1)$, $i=2, \dots, n$ very often  generate
%a subgroup of rank $n-1$ in   the Mordell-Weil group over ${\mathbb Q}$
%of the jacobian of $X_{F,h}/{\mathbb Q}$.
%When such is the case,  One way to prove that many
%choices of $(a_1, \dots, a_n)$ have this property could be as follows.
%Pick first polynomials $a_1(t), \dots, a_n(t)$ in ${\mathbb Q}[t]$ such that
%$\mbox{ord}_t(a_i(t)) = i$. 
\end{remark}

\begin{example} Let $A, B, C \in {\mathbb Z}$ with ${\rm gcd}(A, B, C) = 1$.
Consider the generalized Fermat equation
$$
A x^n + B y^n = Cz^n.
$$
Let $F_{A,B,C,n}/{\mathbb Q}$ denote the projective curve
defined by this equation. 
The reader will easily check that if $F_{A,B,C,n}({\mathbb Q})$ is not empty,
then there exist $A'$, $B'$, and $C'$, such that
$F_{A,B,C,n}$ and $F_{A',B',C',n}$ are isomorphic over ${\mathbb Q}$
and 
 at most one of the coefficients $A'$, $B'$, and $C'$ is
divisible by $p$. We will thus restrict our attention to the case
where $p\nmid AB$. Let ${\mathbb \mu}_n$ denote the group of $n$-th roots of unity in $\overline{\mathbb Q}$. The  curve $F_{A,B,C,n}$  has at least $n^2$ automorphisms:
$$
(x, y) \mapsto (\xi^a_n x, \xi^b_n y), \ \ 1 \le a, b \le n.
$$
In other words, the group ${\mathbb \mu}_n \times {\mathbb \mu}_n$ acts
on $F_{A,B,C,n}/{\mathbb Q}(\xi_n)$.  
Thus, given a single point $P=(x:y:z)$  in $F_{A, B, C,n}({\mathbb Q}(\xi_n))$
 with $xyz \neq 0$, we obtain $n^2$ distinct points in
$F_{A,B,C,n}({\mathbb Q}(\xi_n))$ by considering the orbit of $P$
under ${\mathbb \mu}_n \times {\mathbb \mu}_n$.
When $A = B$ and $P = (x:y:z)$ is such that $x \ne y$, then
we obtain $2n^2$ points in $F_{A,B,C,n}({\mathbb Q}(\xi_n))$ using
the extra automorphism $(x:y:z) \mapsto (y:x:z)$.
\end{example}

\begin{prop}  \label{cor.Fermat}
Let $n = p - 1$ and $p\nmid AB$.  
If the Mordell-Weil rank
over ${\mathbb Q}$ of $F_{A,B,C,p-1}/{\mathbb Q}$ is smaller than $(p-3)/2$, then there exists at most one triple
$(x^{p-1}, y^{p-1}, z^{p-1})$ with $x,y,z \in {\mathbb Z}$,
 $xyz \ne 0$, and such that $Ax^{p-1} + By^{p-1} = Cz^{p-1}$.
Moreover, if $A = B$, then no such triple can exist with $x
\ne y$. 
\end{prop}
\noindent{\bf Proof:} Using \ref{pro.dim}, we find that our hypothesis
implies that the Chabauty rank  of $F_{A,B,C,p-1}$ over ${\mathbb Q}(\xi_{p-1}) $
is less than  the genus of $F_{A,B,C,p-1}$. Thus we may apply \ref{pro.2}
to show that $|F_{A,B,C,p-1}({\mathbb Q}(\xi_{p-1}))| < 2(p-1)^2$.
Proposition \ref{cor.Fermat} follows from the fact that a solution
$P=(x:y:z)$   with $xyz \neq 0$ produces $(p-1)^2$ distinct points in
$F_{A,B,C,p-1}({\mathbb Q}(\xi_{p-1}))$.

\medskip 

Let us say that two solutions $  (x_1, y_1, z_1)$ and $  (x_2,
y_2, z_2)$ of the generalized Fermat equation
are non-equivalent if $(x_1^n, y^n_1, z^n_1)
\ne \lambda(x^n_2, y^n_2, z^n_2)$. 
We shall say that a solution $(x,y,z)$ is non-trivial if $xyz \neq 0$.

Recall that the jacobian of Fermat curve $F_{A,B,C,n}/{\mathbb Q}$ has many quotients.  In particular,
let $a, b \in {\mathbb N}$, and $d:= \gcd(n,a,b)$. Consider the 
curve $Y/{\mathbb Q}$
given by $B^{b/d}y^{n/d} = x^{a/d}(C - Ax)^{b/d}$.  We have a map $F_{A,B,C,n}
\rightarrow Y$, given in coordinates by
$$
(x: y: 1) \mapsto (x^n, x^a y^b).
$$

\begin{cor}  Let $p$ be an odd prime. % and $d:= \gcd(a,b,p-1)$.  
Assume that the curve $A x^{p-1} + By^{p-1} =
Cz^{p-1}$ has two non-trivial non-equivalent solutions $(x_1, y_1, z_1)$ and $(x_2, y_2,
z_2)$ in ${\mathbb Q}^3$ and that  $p \nmid AB$.
Let $a,b \in {\mathbb N}$.
 Let $\ell$ be an odd prime with $\ell \mid
p-1$ and $\ell \nmid ab$.
 Let $Y/{\mathbb Q}$ be given by the equation 
$B^{b} y^{\ell} = x^{a}(C-Ax)^{b}$.
  %and %, when $\ell \mid (p-1)/d$, 
%let $Y'/{\mathbb Q}$ be given by the equation $B^{b/d} y^{\ell} = %x^{a/d}(C-Ax)^{b/d}$.
If $Y$ has positive genus, then its
Mordell-Weil rank over ${\mathbb Q}$ is positive.
\end{cor}

\noindent{\bf Proof:}  Let $F:=F_{A,B,C,p-1}/{\mathbb Q}$.
Our  hypotheses imply that $\Chab(F, {\mathbb Q}(\xi_{p-1}), (p))
= g(F)$.
%the
%jacobian of the curve  has Chabauty rank
%over ${\mathbb Q}(\xi_{p-1})$ at least equal to its genus. 
Thus, any
quotient of ${\rm Jac}(F)/{\mathbb Q}$ defined over ${\mathbb Q}$ has
Chabauty rank over ${\mathbb Q}(\xi_{p-1})$ equal to its
dimension.  Since ${\ell}$ divides $p-1$, the curve $F_{A,B,C,\ell}$
is a surjective image of $F$, the curve $Y/{\mathbb Q}$ is a
surjective image of $F_{A,B,C,\ell}$, and the Jacobian of $Y/{\mathbb
 Q}$ is a quotient of ${\rm Jac}(F)/{\mathbb Q}$.   Hence, the
Chabauty rank of the Jacobian of $Y/{\mathbb
 Q}$ is equal to its dimension and we may apply
\ref{pro.dim} with $s=3$. Indeed, let $f\ell$ denote the smallest positive multiple
of $\ell$ that is greater than or equal to $a+b$. Change $x$ to $x'+m$ for some appropriate
$m$ so that $x'=0$ is not
in the branch locus of the natural map $Y \to {\mathbb P}^1$.
 Then change coordinates to $u:= 1/x'$ and $v = y/{x'}^{f}$
to get an equation for $Y$ of the form $v^{\ell} = u^{f\ell-a-b}(u-\alpha_1)^a(u-\alpha_2)^b $. 
Thus, when $g(Y)>0$, \ref{pro.dim} with $s=3$ shows that 
the Mordell-Weil rank of $Y$  over ${\mathbb Q}$ is greater than or equal to $(s-2)/2 > 0$, which completes our proof.

\medskip
Let $Y'/{\mathbb Q}$ be given by the equation $B^{b/d} y^{\ell} = x^{a/d}(C-Ax)^{b/d}$.
Since  ${\ell}$ divides $(p-1)/d$, the curve $Z$ given by $B^{b/d} y^{n/d} = x^{a/d}(C-Ax)^{b/d}$ is a quotient of $F$, and $Y'$ is a quotient of $Z$.
An isomorphism between the curves $Y'$ and $Y$ is given by $(x,y) \mapsto (x,y^d)$. 
\begin{remark}
It is easy to produce triples $(A, B, C)$ with $ABC \ne 0$
such that $F_{A,B,C,n}({\mathbb Q})$ has 2 non-trivial 
non-equivalent solutions and $p$
divides at most one of $A,B,C$.  Namely, pick appropriate
triples of pairwise coprime integers $(x_1, y_1, z_1)$ and $(x_2,
y_2, z_2)$ with $p \nmid x_1 y_1 z_1 x_2 y_2 z_2$
and
$(x_1^{n}, y^{n}_1, z^{n}_1) \ne (x^{n}_2, y^{n}_2, z^{n}_2)$. 
 Then solve
for $(A, B,C)$ in the simultaneous equations
$$\left\{\begin{array}{c}
Ax_1^{n} + By_1^{n} = Cz_1^{n}\\
Ax_2^{n} + By_2^{n} = Cz^{n}_2.
\end{array}
\right.
$$
  It follows from \ref{cor.Fermat} that for
such a choice of $(A, B, C)$ with $n=p-1$, the Mordell-Weil rank of ${\rm
Jac}(F_{A,B,C,p-1})/{\mathbb Q}$ is greater than or equal to
$(p-3)/2$. On the other hand, it does seem to  follow
from \ref{cor.Fermat} that for
such a choice of $(A, B, C)$ with two solutions $P$ and $Q$,
the ${\mathbb Q}$-rational point $P-Q$ has infinite order in the jacobian
of $F_{A,B,C,p-1}$.

It would be interesting to determine, for any $n$, whether,
when $F_{A,B,C,n}({\mathbb Q})$ contains two non-trivial non-equivalent
points $P$ and $Q$, then $P-Q$ has infinite order in the jacobian of $F_{A,B,C,n}$.
We can construct examples where $P-Q$ has infinite order as follows.
Choose a prime $q>2$ such that
$q \nmid x_1 y_1 z_1 (x_1-y_1)(x_1-z_1)(y_1-z_1)$
and let $(x_2,y_2,z_2) := (x_1, y_1, z_1) + (q,q,q)$.
With these choices, we find that $q \nmid ABC$, where
$(A,B,C)$ is a solution to the above system of equations
with $\gcd(A,B,C) =1$.
It follows that $F_{A,B,C,n}$ has good reduction at $q$.
Moreover, by construction, $P-Q$ is in the kernel of the reduction
mod $q$. Since $q > 2$, the kernel of the reduction does not contain any point
of finite order and, thus, $P-Q$ has infinite order.

When $F_{A,B,C,n}$ has an elliptic quotient $E/{\mathbb Q}$, $P-Q$
has infinite order in the jacobian of $F_{A,B,C,n}$ if the image $P'-Q'$
of $P-Q$ in $E({\mathbb Q})$ has infinite order. For instance,
when $6 \mid n$, we find that $F_{A,B,C,n}$ has  $Ax^3+B = Cz^2$ as quotient,
which can be rewritten as $v^2= u^3 + C^3A^2B$. When $C^3A^2B$ is neither a square, nor
a cube, nor equals $-1$, $E({\mathbb Q})_{tors}$ is trivial. It is easy to check in this case that $P'-Q'$ has infinite order.
\end{remark}

\begin{remark}  Gross and Rohrlich have shown in \cite{G-R}, 2.1,
that if $\ell >7$ is prime, all but three of the $\ell-2$
 Fermat
quotients $y^{\ell} = x^a(1 - x)^b$ have positive
Mordell-Weil ranks over ${\mathbb Q}$.

  Suppose that $A = B$ and $(x:x:z) \in F_{A,A,C,n}({\mathbb Q})$.
Then the curve $F_{A,A,C}$
is isomorphic over ${\mathbb Q}$ to the curve $F_{1,1,2}$,
with the point $(x:x:z)$ corresponding to the point $(1:1:1)$.
The rational points of the curve $x^n + y^n = 2z^n$ are determined
in \cite{D-M}. Other cases are treated in \cite{Sit}.
For conjectures regarding the solutions
of generalized Fermat equations, see \cite{Gra}.

%% Tom
The curve $F_{A,B,C, n}$ is a twist of the curve $F_{1,1,1, n}$.
It is shown in \cite{NewSil}, Thm. 1, that the number of points
over a number field $K$ of a twist $C_{\chi}/K$ of a curve $C/K$
can be bounded in terms of a constant $\gamma(C/K)$ and 
the Mordell-Weil rank of $C_{\chi}/K$.

\end{remark}

\vspace{.25in}
\begin{tabbing}
 
 Department of Mathematics \\
 University of Georgia\\
   Athens, GA 30602, USA\\
 {\tt lorenzini{@}math.uga.edu} \\
{\tt ttucker{@}math.uga.edu} 

\end{tabbing}

\end{document}